\documentclass[11pt,reqno,twoside]{amsart}

\usepackage{amsmath, amssymb, amscd}

\numberwithin{equation}{section}

\theoremstyle{plain}
\newtheorem{thm}{Theorem}
\newtheorem{cor}[thm]{Corollary}

\newtheorem{lem}[thm]{Lemma}
\newtheorem{prop}[thm]{Proposition}
\newtheorem{Def}[thm]{Definition}

\theoremstyle{remark}

\newcounter{sspar}[subsection]
\renewcommand\thesspar{(\thesubsection.\arabic{sspar})}
\newenvironment{sspar}%
    {\par\ \newline
     \vskip-\baselineskip\vskip.1truecm
     \noindent\refstepcounter{sspar}
     \noindent\textbf{\thesspar} \ignorespaces}
    {\vskip-\baselineskip
    \ignorespaces}
\newenvironment{tsspar}%
    {\refstepcounter{sspar}
     \textup{\textbf{\thesspar}} \ignorespaces}
    {\vskip-\baselineskip
    \ignorespaces}
\newcommand{\R}{{\mathbb R}}
\newcommand{\N}{{\mathbb N}}
\newcommand{\C}{{\mathbb C}}
\newcommand{\Z}{{\mathbb Z}}
\newcommand{\Zp}{{\mathbb Z}_{\geq 0}}
\newcommand{\al}{\alpha}

\newcommand{\ga}{\gamma}
\newcommand{\Ga}{\Gamma}
\newcommand{\de}{\delta}
\newcommand{\ep}{\varepsilon}
\newcommand{\si}{\sigma}

\newcommand{\te}{\theta}
\newcommand{\la}{\lambda}
\newcommand{\vp}{\varphi}
\newcommand{\om}{\omega}
\newcommand{\Hi}{\ell^2(\Z)}
\newcommand{\Hp}{\ell^2(\Zp)}
\newcommand{\hf}{{\frac{1}{2}}}
\title{Spectral theory and special functions}
\thanks{Lecture notes for a four hour course in the
SIAM Activity Group ``Orthogonal Polynomials and Special Functions''
Summer School, Laredo, Spain, July 24--28, 2000.  \newline
Last revision: July 5, 2001}
\author{Erik Koelink}

\address{Technische Universiteit Delft, ITS-TWA,
Postbus 5031, 2600 GA Delft, the Netherlands}
\email{koelink@dutiaw4.twi.tudelft.nl}

\begin{document}

\begin{abstract}
A short introduction to the use of the spectral theorem
for self-adjoint operators in the theory of special functions
is given. As the first example,
the spectral theorem is applied to Jacobi operators,
i.e. tridiagonal operators,
on $\Hp$, leading to a proof of Favard's theorem
stating that polynomials
satisfying a three-term recurrence relation are orthogonal
polynomials. We discuss the link to the moment problem.
In the second example, the spectral theorem is
applied to Jacobi operators on $\Hi$. We discuss the
theorem of Masson and Repka linking
the deficiency indices of a Jacobi operator
on $\Hi$ to those of two Jacobi operators on $\Hp$.
For two examples of Jacobi operators on $\Hi$,
namely for the Meixner, respectively Meixner-Pollaczek, functions,
related to the
associated Meixner, respectively Meixner-Pollaczek, polynomials,
and for the second order hypergeometric $q$-difference
operator, we calculate the spectral measure
explicitly. This gives explicit (generalised)
orthogonality relations for hypergeometric
and basic hypergeometric series.
\end{abstract}
\maketitle

\noindent
{\bf Contents.}
\begin{enumerate}
\item[{1.}] Introduction
\item[{2.}] The spectral theorem
\begin{enumerate}
\item[{2.1.}] Hilbert spaces and bounded operators
\item[{2.2.}] The spectral theorem for bounded self-adjoint
operators
\item[{2.3.}] Unbounded self-adjoint operators
\item[{2.4.}] The spectral theorem for unbounded self-adjoint
operators
\end{enumerate}
\item[{3.}] Orthogonal polynomials and Jacobi operators
\begin{enumerate}
\item[{3.1.}] Orthogonal polynomials
\item[{3.2.}] Moment problems
\item[{3.3.}] Jacobi operators
\item[{3.4.}] Unbounded Jacobi operators
\end{enumerate}
\item[{4.}] Doubly infinite Jacobi operators
\begin{enumerate}
\item[{4.1.}] Doubly infinite Jacobi operators
\item[{4.2.}] Relation with Jacobi operators
\item[{4.3.}] The Green kernel
\item[{4.4.}] Example: the Meixner functions
\item[{4.5.}] Example: the basic hypergeometric difference
operator
\end{enumerate}
\item[] References
\end{enumerate}


\newpage
\section{Introduction}

In these lecture notes we give a short introduction to the use
of spectral theory in the theory of special functions. Conversely,
special functions can be used to determine explicitly the spectral
measures of explicit operators on a Hilbert space. The main
ingredient from functional analysis that we are using is the
spectral theorem, both for bounded and unbounded self-adjoint
operators. We recall the main results of the theory in \S 2, hoping
that the lecture notes become more self-contained in this way.
For differential operators this is a very well-known subject,
and one can consult e.g. Dunford and Schwartz \cite{DunfS}.

In \S 3 we discuss Jacobi operators on $\Hp$, and we present
the link to orthogonal polynomials. Jacobi operators on $\Hp$ are
symmetric tridiagonal matrices, and the link to orthogonal
polynomials goes via the three-term recurrence relation
for orthogonal polynomials.
We prove Favard's theorem
in this setting, which is more or less equivalent to the
spectral decomposition of the Jacobi operator involved. We
first discuss the bounded case. Next we discuss the
unbounded case and its link to the classical Hamburger moment
problem. This is very classical, and it can be traced back
to at least Stone's book \cite{Ston}.
This section is much inspired by Akhiezer \cite{Akhi},
Berezanski\u\i\ \cite[\S 7.1]{Bere},
Deift \cite[Ch.~2]{Deif} and Simon \cite{Simo}, and it can
be viewed as an introduction to \cite{Akhi} and \cite{Simo}.
Especially, Simon's paper \cite{Simo} is recommended for
further reading on the subject. See also the recent book
\cite{Tesc} by
Teschl on Jacobi operators and the relation to non-linear
lattices, see also Deift \cite[Ch.~2]{Deif} for the example
of the Toda lattice.
We recall this material on orthogonal polynomials and
Jacobi operators, since it is an important ingredient
in \S 4.

In \S 4 we discuss Jacobi operators on $\Hi$, and we give the
link between a Jacobi operator on $\Hi$ to two Jacobi operators
on $\Hp$ due to Masson and Repka \cite{MassR}, stating
in particular that the Jacobi operator on $\Hi$ is (essentially)
self-adjoint if and only if the two Jacobi operators on $\Hp$
are (essentially) self-adjoint. Next we discuss the example
for the Meixner functions in detail,
following Masson and Repka \cite{MassR}, but the spectral measure
is now completely worked out. In this case the spectral
measure is purely discrete. If we restrict the Jacobi operator
acting on $\Hi$ to a Jacobi operator on $\Hp$, we obtain the
Jacobi operator for the associated Meixner polynomials.
The case of the Meixner-Pollaczek functions is
briefly considered.
As another example we discuss the second order
$q$-hypergeometric difference operator. In this example the
spectral measure has a continuous part and a discrete part.
Here we follow Kakehi \cite{Kake} and \cite[App.~A]{KoelSsu}.
These operators naturally occur in the representation theory
of the Lie algebra $\mathfrak{su}(1,1)$, see \cite{MassR},
or of the quantised universal enveloping algebra
$U_q(\mathfrak{su}(1,1))$, see \cite{KoelSsu}.
Here the action of certain elements from the
Lie algebra or the quantised universal enveloping
algebra is tridiagonal, and one needs to obtain the
spectral resolution. It is precisely this interpretation that
leads to the limit transition discussed in
\ref{4250}.

However, from the point of view of special functions,
the Hilbert space $\Hi$ is generally not the appropriate Hilbert space
to diagonalise a second-order $q$-difference operator $L$.
For the two examples in \S 4 this works nicely, as shown there.
In particular this is true for the ${}_2\vp_1$-series, as
shown in \S 4, see also \cite{CiccKK} for another example.
But for natural extensions of this situation
to higher levels of special functions
this Hilbert space is not good enough. We refer to \cite{KoelSbig} for
the case of the second order $q$-difference operator having
${}_3\vp_2$-series as eigenfunctions corresponding to the big
$q$-Jacobi functions, and to \cite{KoelSAW} for
the case of the second order $q$-difference operator having
${}_8W_7$-series as eigenfunctions corresponding to the Askey-Wilson
functions. For more references and examples of spectral analysis
of second order difference equations we refer to \cite{KoelSNATO}.

There is a huge amount of material on orthogonal polynomials,
and there is a great number of good introductions to orthogonal
polynomials, the moment problem, and the functional analysis
used here. For orthogonal polynomials I have used \cite{Chih},
\cite[Ch.~2]{Deif}
and \cite{VAssc}. For the moment problem there are the
classics by Shohat and Tamarkin \cite{ShohT}
and Stone \cite{Ston}, see also Akhiezer \cite{Akhi},
Simon \cite{Simo} and, of course,
Stieltjes's original paper \cite{Stie}
that triggered the whole subject.
The spectral theorem can be found in many places, e.g.
Dunford and Schwartz \cite{DunfS} and Rudin \cite{Rudi}.

Naturally, there are many more instances of the use of
functional analytic results that can be applied to
special functions. As an example of a more qualitative
question you can wonder how perturbation of the coefficients
in the three-term recurrence relation for orthogonal polynomials
affects the orthogonality measure, i.e. the spectral measure
of the associated Jacobi operator.
Results of this kind can be obtained by using perturbation results
from functional analysis, see e.g. Dombrowski \cite{Domb} and
further references given there.

\noindent
{\it Acknowledgement.} I thank the organisors, especially
Renato \'Alvarez-Nodarse and Francisco Marcell\'an,
of the summer school for inviting me to give lectures
on this subject. Moreover, I thank the participants,
as well as Johan Kustermans,
Hjalmar Rosengren and especially Wolter Groenevelt, for bringing
errors to my attention.

\newpage
\section{The spectral theorem}

\subsection{Hilbert spaces and
bounded operators}

\begin{sspar}\label{205}
A vector space ${\mathcal H}$ over $\C$ is an inner product space if
there exists a mapping $\langle\cdot,\cdot\rangle\colon
{\mathcal H}\times {\mathcal H} \to \C$ such that for all $u,v,w\in{\mathcal H}$
and for all $a,b\in\C$ we have
(i) $\langle a v+b w, u\rangle = a \langle v,u\rangle
+b \langle w,u\rangle$, (ii) $\langle u,v\rangle = \overline{
\langle v,u\rangle}$, and (iii) $\langle v,v\rangle \geq 0$ and
$\langle v,v\rangle =0$ if and only if $v=0$. With
the inner product we associate the norm $\| v\|=\|v\|_{\mathcal H}
=\sqrt{\langle v,v \rangle}$, and the topology from the
corresponding metric $d(u,v)=\|u-v\|$. The standard inequality
is the Cauchy-Schwarz inequality;
$|\langle u,v\rangle| \leq \| u\| \|v\|$.
A Hilbert space ${\mathcal H}$
is a complete inner product space, i.e. for
any Cauchy sequence $\{x_n\}_n$ in ${\mathcal H}$, i.e.
$\forall \ep>0$ $\exists N\in\N$ such that for
all $n,m\geq N$ $\| x_n-x_m\|<\ep$, there exists an element
$x\in{\mathcal H}$ such that $x_n$ converges to $x$.
In these notes all Hilbert spaces are separable, i.e. there
exists a denumerable set of basis vectors.
\end{sspar}

\begin{sspar}\label{210}
{\it Example.} $\Hi$, the space of square summable
sequences $\{ a_k\}_{k\in\Z}$, and $\Hp$, the space of square
summable sequences $\{ a_k\}_{k\in\Zp}$, are Hilbert spaces.
The inner product is given by
$\langle \{ a_k\}, \{ b_k\} \rangle = \sum_k a_k\overline{b_k}$.
An orthonormal basis is given by the sequences $e_k$ defined
by $(e_k)_l=\de_{k,l}$, so we identify $\{a_k\}$ with
$\sum_k a_ke_k$.
\end{sspar}

\begin{sspar}\label{220} {\it Example.}
We consider a positive Borel measure $\mu$ on the real line
$\R$ such that all moments exist, i.e.
$\int_\R |x|^m\, d\mu(x) < \infty$ for all $m\in\Zp$. Without loss
of generality we assume that $\mu$ is a probability measure,
$\int_\R d\mu(x)=1$.
By $L^2(\mu)$ we denote the space of square integrable functions
on $\R$, i.e. $\int_\R |f(x)|^2\, d\mu(x)<\infty$.
Then $L^2(\mu)$ is a Hilbert space (after identifying
two functions $f$ and $g$ for which
$\int_\R |f(x)-g(x)|^2\, d\mu(x)=0$) with respect to the
inner product $\langle f,g\rangle =\int_\R f(x)\overline{g(x)}\,
d\mu(x)$.
In case $\mu$ is a finite sum of discrete Dirac measures,
we find that $L^2(\mu)$ is finite dimensional.
\end{sspar}

\begin{sspar}\label{230} An operator $T$ from a Hilbert space ${\mathcal H}$ into
another Hilbert space ${\mathcal K}$ is linear if
for all $u,v\in {\mathcal H}$ and for all $a,b\in\C$ we have
$T(au+bv)=aT(u)+bT(v)$. An operator $T$ is bounded if there exists
a constant $M$ such that $\| Tu\|_{\mathcal K}\leq M\|u\|_{\mathcal H}$
for all $u\in {\mathcal H}$. The smallest $M$ for which this holds is
the norm, denoted by $\|T\|$, of $T$. A bounded linear operator
is continuous. The adjoint of a bounded linear operator $T\colon
{\mathcal H}\to{\mathcal K}$ is a map $T^\ast\colon {\mathcal K}\to {\mathcal H}$
with $\langle Tu,v\rangle_{\mathcal K} =
\langle u, T^\ast v\rangle_{\mathcal H}$. We call $T\colon
{\mathcal H}\to {\mathcal H}$ self-adjoint if $T^\ast=T$.
It is unitary if $T^\ast T = {\mathbf 1}_{\mathcal H}$ and
$TT^\ast = {\mathbf 1}_{\mathcal K}$.
A projection $P\colon {\mathcal H}\to {\mathcal H}$ is a linear map
such that $P^2=P$.
\end{sspar}

\begin{sspar}\label{240} We are also
interested in unbounded linear operators.
In that case we denote $(T, {\mathcal D}(T))$, where ${\mathcal D}(T)$, the
domain of $T$, is
a linear subspace of ${\mathcal H}$ and $T\colon {\mathcal D}(T)\to {\mathcal H}$.
Then $T$ is densely defined if the closure of ${\mathcal D}(T)$ equals
${\mathcal H}$. All unbounded operators that we consider in these
notes are densely defined. If the operator $(T-z)$, $z\in\C$, has an
inverse $R(z)=(T-z)^{-1}$ which is densely defined and is bounded,
so that $R(z)$,
the resolvent operator,
extends to a bounded linear operator on ${\mathcal H}$,
then we call $z$ a regular value. The set of all regular values
is the resolvent set $\rho(T)$. The complement of the resolvent
set $\rho(T)$ in $\C$ is the spectrum $\si(T)$ of $T$.
The point spectrum is the subset of the spectrum for which
$T-z$ is not one-to-one. In this case there exists a vector
$v\in{\mathcal H}$ such that $(T-z)v=0$, and $z$ is an eigenvalue.
The continuous spectrum consists of the points $z\in\si(T)$ for
which $T-z$ is one-to-one, but for which $(T-z){\mathcal H}$ is dense
in ${\mathcal H}$,
but not equal to ${\mathcal H}$. The remaining part of the spectrum
is the residual spectrum. For self-adjoint operators, both
bounded and unbounded, see \ref{2110},
the spectrum only consists of the discrete and continuous spectrum.
\end{sspar}

\begin{sspar}\label{245} For a bounded operator $T$ the spectrum $\si(T)$
is a compact subset of the disk of radius $\| T\|$. Moreover,
if $T$ is self-adjoint, then $\si(T)\subset \R$, so that
$\si(T)\subset [-\|T\|, \|T\|]$ and the spectrum consists
of the point spectrum and the continuous spectrum.
\end{sspar}

\newpage
\subsection{The spectral theorem
for bounded self-adjoint operators}

\begin{sspar}\label{250} A resolution of the identity, say $E$, of a Hilbert space
${\mathcal H}$ is a projection valued Borel measure on $\R$ such that
for all Borel sets $A,B\subseteq \R$ we have
(i) $E(A)$ is a self-adjoint projection, (ii) $E(A\cap B)=
E(A)E(B)$, (iii) $E(\emptyset)=0$,
$E(\R)={\mathbf 1}_{\mathcal H}$,
(iv) $A\cap B=\emptyset$ implies $E(A\cup B)=E(A)+E(B)$,
and (v) for all $u,v\in{\mathcal H}$ the map $A\mapsto
E_{u,v}(A)=\langle E(A)u,v\rangle$ is a complex Borel measure.
\end{sspar}

\begin{sspar}\label{260} A generalisation of the spectral
theorem for matrices
is the following theorem for bounded self-adjoint operators,
see \cite[\S X.2]{DunfS}, \cite[\S 12.22]{Rudi}.

\begin{thm}\nonumber{\rm (Spectral theorem)} Let $T\colon
{\mathcal H}\to {\mathcal H}$ be a bounded self-adjoint linear map, then
there exists a unique resolution of the identity such that
$T=\int_\R t \, dE(t)$, i.e. $\langle Tu,v\rangle =\int_\R t \,
dE_{u,v}(t)$. Moreover, $E$ is supported on the
spectrum $\si(T)$, which is contained in the interval
$[-\|T\|, \|T\|]$.
Moreover, any of the spectral projections $E(A)$, $A\subset \R$
a Borel set, commutes with $T$.
\end{thm}

A more general theorem of this kind holds for normal
operators, i.e. for those operators satisfying $T^\ast T=TT^\ast$.
\end{sspar}

\begin{sspar}{\label{270}} Using the spectral theorem we define for
any continuous
function $f$ on the spectrum $\si(T)$ the operator $f(T)$ by
$f(T)=\int_\R f(t) \, dE(t)$, i.e.
$\langle f(T)u,v\rangle =\int_\R f(t) \, dE_{u,v}(t)$. Then $f(T)$
is bounded operator with norm equal to
the supremum norm of $f$ on the spectrum of $T$, i.e.
$\| f(T)\| = \sup_{x\in\si(T)} |f(x)|$.
This is known as the functional calculus for self-adjoint operators.
In particular, for
$z\in\rho(T)$ we see that $f\colon x\mapsto (x-z)^{-1}$
is continuous on the spectrum, and the corresponding operator is
just the resolvent operator $R(z)$ as in \ref{240}. The
functional calculus can be extended to measurable functions,
but then $\| f(T)\| \leq \sup_{x\in\si(T)} |f(x)|$.
\end{sspar}

\begin{sspar}{\label{280}} The spectral measure can be obtained from the
resolvent operators by the Stiel\-tjes-Perron inversion formula,
see \cite[Thm. X.6.1]{DunfS}.

\begin{thm} The spectral measure of the open interval
$(a,b)\subset \R$ is given by
$$
E_{u,v}\bigl( (a,b)\bigr) = \lim_{\de\downarrow 0}
\lim_{\ep\downarrow 0} \frac{1}{2\pi i}
\int_{a+\de}^{b-\de}
\langle R(x+i\ep)u,v\rangle - \langle R(x-i\ep)u,v\rangle \, dx.
$$
The limit holds in the strong operator topology, i.e.
$T_nx\to Tx$ for all $x\in{\mathcal H}$.
\end{thm}
\end{sspar}

\subsection{Unbounded self-adjoint operators}

\begin{sspar}{\label{290}} Let $(T, {\mathcal D}(T))$, with ${\mathcal D}(T)$ the
domain of $T$, be a densely defined unbounded operator on
${\mathcal H}$, see \ref{240}. We can now define the
adjoint operator $(T^\ast, {\mathcal D}(T^\ast))$ as follows.
First define
$$
{\mathcal D}(T^\ast) = \{ v\in {\mathcal H}\mid
u\mapsto \langle Tu,v\rangle \text{\ is continuous on
${\mathcal D}(T)$} \}.
$$
By the density of ${\mathcal D}(T)$
the map $u\mapsto \langle Tu,v\rangle$
for $v\in {\mathcal D}(T^\ast)$ extends to a continuous linear
functional $\om\colon
{\mathcal H}\to \C$, and by the Riesz representation
theorem there exists a unique $w\in{\mathcal H}$ such that
$\om(u)=\langle u,w\rangle$ for all $u\in {\mathcal H}$.
Now the adjoint $T^\ast$ is defined by $T^\ast v =w$, so that
$$
\langle Tu,v\rangle = \langle u, T^\ast v\rangle \qquad
\forall\, u\in{\mathcal D}(T),\, \forall \, v\in {\mathcal D}(T^\ast).
$$
\end{sspar}

\begin{sspar}{\label{2100}} If $T$ and $S$ are unbounded operators
on ${\mathcal H}$, then $T$ extends $S$, notation
$S\subset T$, if ${\mathcal D}(S)\subset {\mathcal D}(T)$ and
$Sv=Tv$ for all $v\in{\mathcal D}(S)$. Two unbounded operators
$S$ and $T$ are equal, $S=T$, if $S\subset T$ and $T\subset S$,
or $S$ and $T$ have the same domain and act in the same way.
In terms of the graph
$$
{\mathcal G}(T) = \{ (u, Tu)\mid u\in {\mathcal D}(T)\} \subset
{\mathcal H}\times {\mathcal H}
$$
we see that $S\subset T$ if and only if ${\mathcal G}(S) \subset
{\mathcal G}(T)$. An operator $T$ is closed if its graph is
closed in the product topology of ${\mathcal H}\times {\mathcal H}$.
The adjoint of a densely defined operator is a closed
operator, since the graph of the adjoint is given as
$$
{\mathcal G}(T^\ast) = \{ (-u,Tu)\mid u\in {\mathcal D}(T)\}^\perp,
$$
for the inner product
$\langle (u,v),(x,y)\rangle=\langle u,x\rangle
+\langle v,y\rangle$ on ${\mathcal H}\times {\mathcal H}$, see
\cite[13.8]{Rudi}.
\end{sspar}

\begin{sspar}{\label{2110}}
A densely defined operator is symmetric if
$T\subset T^\ast$, or, using the definition in \ref{290},
$$
\langle Tu,v\rangle = \langle u,Tv\rangle , \qquad\forall \
u,v\in{\mathcal D}(T).
$$
A densely defined operator is self-adjoint if $T=T^\ast$,
so that a self-adjoint operator is closed. The spectrum of
an unbounded self-adjoint operator is contained in $\R$.
Note that ${\mathcal D}(T) \subset {\mathcal D}(T^\ast)$, so that
${\mathcal D}(T^\ast)$ is a dense subspace and taking
the adjoint once more gives
$(T^{\ast\ast}, {\mathcal D}(T^{\ast\ast}))$
as the minimal closed extension
of $(T, {\mathcal D}(T))$, i.e. any densely defined symmetric operator
has a closed extension. We have
$T\subset T^{\ast\ast}\subset T^\ast$.
We say that the densely defined symmetric operator is
essentially self-adjoint if its closure is self-adjoint, i.e.
if $T\subset T^{\ast\ast} = T^\ast$.
\end{sspar}

\begin{sspar}{\label{2120}}
In general, a densely defined symmetric operator $T$
might not have self-adjoint extensions. This can be measured
by the deficiency indices. Define for $z\in \C\backslash\R$
the eigenspace
$$
N_z = \{ v\in {\mathcal D}(T^\ast)\mid T^\ast v=z\, v\}.
$$
Then $\dim N_z$ is constant for $\Im z>0$ and for $\Im z<0$,
\cite[Thm. XII.4.19]{DunfS}, and we put $n_+=\dim N_i$
and $n_-=\dim N_{-i}$. The pair $(n_+,n_-)$ are
the deficiency indices for the densely defined symmetric
operator $T$. Note that if $T$ commutes with complex conjugation
of the Hilbert space ${\mathcal H}$ then we automatically have
$n_+=n_-$. Note furthermore that if $T$ is self-adjoint then
$n_+=n_-=0$, since a self-adjoint operator cannot have
non-real eigenvalues. Now the following
holds, see \cite[\S XII.4]{DunfS}.

\begin{prop} Let $(T, {\mathcal D}(T))$
be a densely defined symmetric operator. \par\noindent
{\rm (i)}
${\mathcal D}(T^\ast)= {\mathcal D}(T^{\ast\ast})
\oplus N_i \oplus N_{-i}$,
as an orthogonal direct sum with respect to the graph norm for
$T^\ast$ from
$\langle u,v\rangle_{T^\ast} = \langle u,v\rangle +
\langle T^\ast u,T^\ast v\rangle$. As a direct sum,
${\mathcal D}(T^\ast)= {\mathcal D}(T^{\ast\ast}) + N_z + N_{\bar z}$
for general $z\in\C\backslash\R$.
\par\noindent
{\rm (ii)} Let $U$ be an isometric bijection $U\colon N_i\to N_{-i}$
and define $(S, {\mathcal D}(S))$ by
$$
{\mathcal D}(S) = \{ u + v +Uv\mid u\in {\mathcal D}(T^{\ast\ast}), \
v\in N_i\}, \quad Sw=T^\ast w,
$$
then $(S, {\mathcal D}(S))$ is a self-adjoint extension
of $(T,{\mathcal D}(T))$,
and every self-adjoint extension of $T$ arises in this way.
\end{prop}

In particular, $T$ has self-adjoint extensions if and only if
the deficiency indices are equal; $n_+=n_-$.
$T^{\ast\ast}$ is a closed symmetric extension of
$T$. We can also
characterise the domains of the self-adjoint extensions of $T$
using the sesquilinear form
$$
B(u,v) =  \langle T^\ast u, v\rangle -
\langle u, T^\ast v\rangle, \qquad u,v\in {\mathcal D}(T^\ast),
$$
then ${\mathcal D}(S)=\{ u\in {\mathcal D}(T^\ast)\mid B(u,v)=0
, \ \forall v\in {\mathcal D}(S)\}$.
\end{sspar}

\subsection{The spectral theorem
for unbounded self-adjoint operators}

\begin{sspar}{\label{2130}} With all the preparations of the previous
subsection the Spectral
Theorem \ref{260} goes through in the unbounded setting,
see \cite[\S XII.4]{DunfS}, \cite[Ch.~13]{Rudi}.

\begin{thm}{\rm (Spectral theorem)} Let $T\colon
{\mathcal D}(T)\to {\mathcal H}$ be an unbounded self-adjoint linear map,
then there exists a unique resolution of the identity such that
$T=\int_\R t \, dE(t)$, i.e. $\langle Tu,v\rangle =\int_\R t \,
dE_{u,v}(t)$ for $u\in{\mathcal D}(T)$, $v\in {\mathcal H}$.
Moreover, $E$ is supported on the
spectrum $\si(T)$, which is contained in $\R$.
For any bounded operator $S$ that satisfies $ST\subset TS$ we
have $E(A)S=SE(A)$, $A\subset \R$
a Borel set.
Moreover, the Stieltjes-Perron inversion formula
\textup{\ref{280}} remains valid;
$$
E_{u,v}\bigl( (a,b)\bigr) = \lim_{\de\downarrow 0}
\lim_{\ep\downarrow 0} \frac{1}{2\pi i}
\int_{a+\de}^{b-\de}
\langle R(x+i\ep)u,v\rangle - \langle R(x-i\ep)u,v\rangle \, dx.
$$
\end{thm}
\end{sspar}

\begin{sspar}{\label{2140}} As in \ref{270} we can now define
$f(T)$ for any measurable function $f$ by
$$
\langle f(T)u,v\rangle = \int_\R f(t)\, dE_{u,v}(t),
\qquad u\in {\mathcal D}(f(T)), \ v\in {\mathcal H},
$$
where ${\mathcal D}(f(T)) = \{ u\in {\mathcal H}\mid
\int_\R |f(t)|^2\, dE_{u,u}(t) <\infty\}$ is the domain
of $f(T)$. This makes $f(T)$ into a densely defined closed
operator. In particular, if $f\in L^\infty(\R)$, then
$f(T)$ is a continuous operator, by the closed graph
theorem. This in particular applies to $f(x)=(x-z)^{-1}$,
$z\in\rho(T)$, which gives the resolvent operator.
\end{sspar}

\newpage
\section{Orthogonal polynomials and Jacobi operators}

\subsection{Orthogonal polynomials}

\begin{sspar}{\label{310}}
Consider the Hilbert space $L^2(\mu)$ as in Example~\ref{220}.
Assume that all moments exist,
so that all polynomials are integrable. In applying the
Gram-Schmidt orthogonalisation process to the sequence
$\{ 1,x,x^2,x^3,\ldots\}$ we may end up in one of the
following situations: (a) the polynomials are linearly
dependent in $L^2(\mu)$, or (b) the polynomials are
linearly independent in $L^2(\mu)$.
In case (a) it follows that there is a non-zero
polynomial $p$ such that $\int_\R |p(x)|^2\, d\mu(x)=0$. This
implies that $\mu$ is a finite sum of Dirac measures at the zeros
of $p$. From now on we exclude this case, but the reader may
consider this case him/herself. In case (b) we end up with a
set of orthonormal polynomials as in the following definition.

\begin{Def} A sequence of polynomials
$\{p_n\}_{n=0}^\infty$ with $\deg(p_n)=n$
is a set of orthonormal polynomials
with respect to $\mu$ if $\int_\R p_n(x) p_m(x) \, d\mu(x) =
\de_{n,m}$.
\end{Def}

Note that the polynomials $p_n$ are real-valued for $x\in\R$, so
that its coefficients are real. Moreover, from the Gram-Schmidt
process it follows that the leading coefficient is positive.
\end{sspar}

\begin{sspar}{\label{350}}
Note that only the moments $m_k=\int_\R x^k\, d\mu(x)$
of $\mu$ play a role in the orthogonalisation process.
The Stieltjes transform of the measure $\mu$
defined by $w(z) = \int_\R (x-z)^{-1}\, d\mu(x)$,
$z\in\C\backslash\R$, can be considered as a generating
function for the moments of $\mu$. Indeed, formally
\begin{equation}
w(z) = \frac{-1}{z} \int_\R \frac{1}{1-x/z}\, d\mu(x) =
\frac{-1}{z} \sum_{k=0}^\infty \int_\R \bigl( \frac{x}{z}\bigr)^k
\, d\mu(x) = -\sum_{k=0}^\infty \frac{m_k}{z^{k+1}}.
\label{eq315}
\end{equation}
In case $\text{supp}(\mu)\subseteq [-A,A]$ we see that
$|m_k|\leq 2A^k$ implying that the series in (\ref{eq315})
is absolutely convergent for $|z|>A$. In this case we see that the
Stieltjes transform $w(z)$ of $\mu$ is completely determined by
the moments of $\mu$.
In general, this expansion has to be interpreted as an asymptotic
expansion of the Stieltjes transform $w(z)$ as $|z|\to\infty$.
We now give a proof of the Stieltjes inversion formula,
cf. \ref{280}.

\begin{prop} Let $\mu$ be a probability measure
with finite moments, and let $w(z) = \int_\R (x-z)^{-1}\, d\mu(x)$
be its Stieltjes transform, then
$$
\lim_{\ep\downarrow 0} \frac{1}{\pi} \int_a^b
\Im \bigl(w(x+i\ep)\bigr)\, dx =  \mu\bigl( (a,b)\bigr)
+\hf \mu(\{ a\}) + \hf \mu(\{b\}).
$$
\end{prop}

\begin{proof} Observe that
\begin{equation*}
\begin{split}
2i \Im \bigl( w(z)\bigr) &=
w(z)-\overline{w(z)} = w(z)-w(\bar z) =
\int_\R \frac{1}{x-z} - \frac{1}{x-\bar z}\, d\mu(x)  \\
&= \int_\R \frac{z-\bar z}{|x-z|^2}\, d\mu(x) = 2i
\int_\R \frac{\Im z}{|x-z|^2}\, d\mu(x),
\end{split}
\end{equation*}
so that
$$
\Im \bigl( w(x+i\ep)\bigr)
= \int_\R \frac{\ep}{|s-(x+i\ep)|^2}\, d\mu(s)
= \int_\R \frac{\ep}{(s-x)^2+\ep^2}\, d\mu(s).
$$
Integrating this expression and interchanging integration,
which is allowed since the integrand is positive, gives
\begin{equation}
\int_a^b  \Im \bigl(w(x+i\ep)\bigr)\, dx =
\int_\R \int_a^b \frac{\ep}{(s-x)^2+\ep^2}\, dx\, d\mu(s).
\label{eq317}
\end{equation}
The inner integration can be carried out easily;
$$
\chi_\ep(s) = \int_a^b \frac{\ep}{(s-x)^2+\ep^2}\, dx =
\int_{(a-s)/\ep}^{(b-s)/\ep} \frac{1}{1+y^2}\, dy =
\arctan y \Big\vert_{(a-s)/\ep}^{(b-s)/\ep}
$$
by $y=(x-s)/\ep$. It follows that $0\leq \chi_\ep(s)\leq \pi$,
and
$$
\lim_{\ep\downarrow 0}\chi_\ep(s) =
\begin{cases}\pi, & \text{for $a<s<b$,} \\
\hf \pi, & \text{for $s=a$ or $s=b$.}\end{cases}
$$
It suffices to show that we can interchange integration and
the limit $\ep\downarrow 0$ in (\ref{eq317}). This
follows from Lebesgue's dominated convergence theorem since
$\mu$ is a probability measure and $0\leq \chi_\ep(s)\leq \pi$.
\end{proof}

As a corollary to the proof we get, cf. \ref{280}, \ref{2130},
$$
\lim_{\ep\downarrow 0} \lim_{\de\downarrow 0}
\frac{1}{\pi} \int_{a+\de}^{b-\de}
\Im \bigl(w(x+i\ep)\bigr)\, dx =  \mu\bigl( (a,b)\bigr).
$$

We need the following extension of this inversion formula
in \ref{3380}. For a polynomial $p$ with real coefficients
we have
\begin{equation}
\lim_{\ep\downarrow 0}
\frac{1}{\pi} \int_a^b
\Im \bigl(p(x+i\ep)w(x+i\ep)\bigr)\, dx =  \int_{(a,b)} p(x)\,d\mu(x)
+\frac{1}{2} p(a)\mu (\{a\}) +\frac{1}{2} p(b)\mu (\{b\}).
\label{eqextra1}
\end{equation}
We indicate how the proof of the proposition can be extended
to obtain (\ref{eqextra1}).
Start with
$$
\Im\bigl(p(x+i\ep)w(x+i\ep)\bigr) = \int_\R
\frac{(s-x)\Im\bigl(p(x+i\ep)\bigr) +\ep \Re\bigl(
p(x-i\ep)\bigr)}{(s-x)^2+\ep^2}\, d\mu(s).
$$
Integrate this expression with respect to $x$
and interchange summations, which is justified
since $\Im\bigl(p(x+i\ep)\bigr)$
and $\Re\bigl(p(x-i\ep)\bigr)$ are bounded on
$(a,b)$. This time we have to evaluate two integrals.
The first integral
$$
\int_a^b
\frac{(s-x)\Im\bigl(p(x+i\ep)\bigr)}{(s-x)^2+\ep^2}\, dx
= \int_{(a-s)/\ep}^{(b-s)/\ep} \frac{-y\, \Im\bigl(p(\ep
y+s+i\ep)\bigr)}{y^2+1}\, dy
$$
can be estimated,
using $\Im\bigl(p(x+i\ep)\bigr)={\mathcal O}(\ep)$
uniformly on $[a,b]$, by
$$
\ep\, M\, \int_{(a-s)/\ep}^{(b-s)/\ep} \frac{y}{y^2+1}\, dy
= \ep M \ln \sqrt{y^2+1} \Big\vert_{(a-s)/\ep}^{(b-s)/\ep}.
$$
This term tends to zero independently of $a$, $b$ and $s$,
since $\ep\ln(\frac{A^2}{\ep^2}+1) = -2\ep\ln\ep +
\ep\ln(A^2+\ep^2)$ which can be estimated by
${\mathcal O}(\ep\ln\ep)$ with a constant independent
of $A$.
The other integral can be dealt with as in the proof
of the proposition. Next Lebesgue's dominated
convergence theorem can be applied and (\ref{eqextra1})
follows.
\end{sspar}

\begin{sspar}{\label{380}} The following theorem describes
the fundamental property of orthogonal polynomials in these notes.

\begin{thm} {\rm (Three term recurrence
relation)} Let $\{ p_k\}_{k=0}^\infty$
be a set of orthonormal polynomials in $L^2(\mu)$,
then there exist sequences $\{ a_k\}_{k=0}^\infty$,
$\{ b_k\}_{k=0}^\infty$,  with $a_k>0$ and $b_k\in\R$, such
that
\begin{eqnarray}
x\, p_k(x) &=& a_kp_{k+1}(x) + b_kp_k(x) + a_{k-1} p_{k-1}(x), \qquad
k\geq 1, \label{eq330}\\
x\, p_0(x) &=& a_0p_1(x) + b_0p_0(x). \label{eq340}
\end{eqnarray}
Moreover, if $\mu$ is compactly supported, then the
coefficients $a_k$ and $b_k$ are bounded.
\end{thm}

Note that (\ref{eq330}), (\ref{eq340}) together with the initial
condition $p_0(x)=1$ completely determine the polynomials $p_k(x)$
for all $k\in \N$.

\begin{proof} The degree of $xp_k(x)$ is $k+1$, so there exist
constants $c_i$ such that
$x\, p_k(x) = \sum_{i=0}^{k+1} c_i\, p_i(x)$. By the orthonormality
properties of $p_k$ it follows that
$$
c_i = \int_\R p_i(x)xp_k(x)\, d\mu(x).
$$
Since the degree
of $xp_i(x)$ is $i+1$, we see that $c_i=0$ for $i+1<k$. Then
$$
b_k=c_k=\int_\R x \bigl(p_k(x))\bigr)^2\, d\mu(x) \in\R.
$$
Moreover,
$c_{k+1} = \int_\R p_{k+1}(x)xp_k(x)\, d\mu(x)$ and
$c_{k-1} = \int_\R p_{k-1}(x)xp_k(x)\, d\mu(x)$ display the
required structure for the other coefficients.
The positivity of $a_k$ follows by considering the leading
coefficient.

For the last statement we observe that
\begin{equation*}
\begin{split}
|a_k| &=\big\vert \int_\R xp_{k+1}(x)p_k(x)\, d\mu(x)\big\vert
\leq \int_\R |p_{k+1}(x)| |p_k(x)|\, d\mu(x) \,
\sup_{x\in\text{supp}(\mu)} |x| \\ &\leq
\|p_{k+1}\|_{L^2(\mu)} \|p_k\|_{L^2(\mu)}
\sup_{x\in\text{supp}(\mu)} |x| =
\sup_{x\in\text{supp}(\mu)} |x|<\infty,
\end{split}
\end{equation*}
since $\|p_k\|_{L^2(\mu)}=1$ and $\text{supp}(\mu)$ is compact.
In the second inequality we have used the Cauchy-Schwarz
inequality \ref{205}. Similarly,
$$
|b_k| \leq \|p_k\|_{L^2(\mu)}^2
\sup_{x\in\text{supp}(\mu)} |x| =
\sup_{x\in\text{supp}(\mu)} |x|<\infty
$$
gives the estimate on the coefficients $b_k$.
\end{proof}
\end{sspar}

\begin{sspar}{\label{3100}}
We observed that (\ref{eq330}) and (\ref{eq340}) together with an
initial condition for the degree zero component completely
determines a solution for the recurrence (\ref{eq330}), (\ref{eq340}).
We can also generate solutions of (\ref{eq330}) by specifying the
initial values for $k=0$ and $k=1$.
{}From now on we let $r_k(x)$ be the sequence of
polynomials generated by
(\ref{eq330}) subject to the initial conditions $r_0(x)=0$ and
$r_1(x)=a_0^{-1}$. Then $r_k$ is a polynomial of degree $k-1$,
and (\ref{eq340}) is not valid. The polynomials
$\{ r_k\}_{k=0}^\infty$ are the associated polynomials.

\begin{lem} The associated polynomial $r_k$ is given by
$$
r_k(x) = \int_\R \frac{p_k(x)-p_k(y)}{x-y} \, d\mu(y).
$$
\end{lem}

\begin{proof} It suffices to show that the right hand side,
denoted temporarily by $q_k(x)$,
satisfies the recurrence (\ref{eq330}) together with
the initial conditions. Using Theorem \ref{380} for $p_k(x)$
and the definition of $q_k(x)$ we obtain
\begin{equation*}
\begin{split}
x\, q_k(x) =& a_kq_{k+1}(x) + b_kq_k(x) + a_{k-1} q_{k-1}
 \\ &+ \int_\R \frac{a_kp_{k+1}(y) + b_kp_k(y) + a_{k-1}
p_{k-1}(y)-xp_k(y)}{x-y}\, d\mu(y).
\end{split}
\end{equation*}
Using Theorem \ref{380} again shows that the integral
equals  $-\int_\R p_k(y)\, d\mu(y)$, which is zero for $k\geq 1$
and $-1$ for $k=0$ by the orthogonality properties. Hence,
(\ref{eq330}) is satisfied. Using $p_0(x)=1$ we find $q_0(x)=0$ and
using $p_1(x) = a_0^{-1}(x-b_0)$ gives $q_1(x)=a_0^{-1}$.
\end{proof}

Considering
$$
\int_\R \frac{p_k(x)}{x-z}\, d\mu(x) =
\int_\R \frac{p_k(x)-p_k(z)}{x-z} \, d\mu(x) +
p_k(z) \int_\R \frac{1}{x-z} \, d\mu(x)
$$
immediately proves the following corollary.

\begin{cor} Let $z\in\C\backslash\R$ be fixed.
The $k$-th coefficient
with respect to the orthonormal set $\{p_k\}_{k=0}^\infty$
in $L^2(\mu)$ of $x\mapsto (x-z)^{-1}$ is given by
$w(z)p_k(z)+r_k(z)$. Hence,
$$
\sum_{k=0}^\infty |w(z)p_k(z)+r_k(z)|^2 \leq \int_\R |x-z|^{-2}\,
d\mu(x)<\infty.
$$
\end{cor}

The inequality follows from the Bessel inequality. If the
$\{p_n\}_{n=0}^\infty$ is an orthonormal basis of $L^2(\mu)$
then we can write an equality by Parseval's identity.
\end{sspar}

\begin{sspar}{\label{3150}}
Since $r_k(y)$ is another solution to (\ref{eq330}),
multiplying (\ref{eq330}) by $r_k(y)$ and (\ref{eq330}) for $r_k(y)$
by $p_k(x)$, subtracting leads to
\begin{equation}
\begin{split}
(x-y)p_k(x)r_k(y) =&
a_k\bigl( p_{k+1}(x)r_k(y) -p_k(x)r_{k+1}(y)\bigr)
\\ &-  a_{k-1}\bigl( p_k(x)r_{k-1}(y) -p_{k-1}(x)r_k(y)\bigr)
\end{split}
\label{eq349}
\end{equation}
for $k\geq 1$.
Taking $x=y$ in (\ref{eq349}), we see that the Wronskian,
or Casorati determinant,
\begin{equation}
[p,r]_k(x)= a_k\bigl( p_{k+1}(x)r_k(x) -p_k(x)r_{k+1}(x)\bigr)
\label{eq351}
\end{equation}
is independent of $k\in\Zp$, and taking $k=0$ gives
$[p,r]_k(x)=[p,r]=-1$.
This also shows that $p_k$ and $r_k$ are linearly
independent solutions to (\ref{eq330}).

On the other hand, replacing $r_k$ by $p_k$ and summing we get
the Christoffel-Darboux formula
\begin{equation}
(x-y) \sum_{k=0}^{n-1} p_k(x)p_k(y) = a_{n-1}\bigl(
p_n(x)p_{n-1}(y) - p_{n-1}(x)p_n(y)\bigr).
\label{eq360}
\end{equation}
The case $x=y$ is obtained after dividing by $x-y$ and in the
right hand side letting $y\to x$. This gives
$$
\sum_{k=0}^{n-1} p_k(x)^2  = a_{n-1}\bigl(
p_n'(x)p_{n-1}(x) - p_{n-1}'(x)p_n(x)\bigr).
$$
\end{sspar}

\subsection{Moment problems}

\begin{sspar}{\label{3250}}
The moment problem consists of the following two
questions:

\begin{enumerate}
\item Given a sequence $\{m_0,m_1, m_2, \ldots\}$, does there
exist a positive Borel measure $\mu$ on $\R$
such that $m_k=\int x^k\, d\mu(x)$?
\item If the answer to problem 1 is yes, is the measure obtained
unique?
\end{enumerate}

\noindent
Note that we can assume without loss of generality that $m_0=1$.
This is always assumed.
\end{sspar}

\begin{sspar}{\label{3260}}
In case $\text{\rm supp}(\mu)$ is required to be
contained in a finite interval we
speak of the Haussdorf moment problem (1920). In case
$\text{\rm supp}(\mu)$ is to be contained in $[0,\infty)$ we
speak of the Stieltjes moment problem (1894). Finally, without
a condition on the support, we speak of the Hamburger moment
problem (1922).
Here, a moment problem is always a Hamburger moment problem.
The answer to question 1 can be given completely in terms of
positivity requirements of matrices composed of the $m_i$'s,
see Akhiezer \cite{Akhi}, Shohat and Tamarkin \cite{ShohT},
Simon \cite{Simo}, Stieltjes \cite{Stie}, Stone \cite{Ston}.
In these notes we only discuss an answer to question 2,
see \S 3.4.

In case the answer to question 2 is affirmative, we speak of
a determinate moment problem and otherwise of an
indeterminate moment problem. So for an indeterminate moment
problem we have a convex set of probability measures on $\R$
solving the same moment problem. The fact that this may happen
has been observed first by Stieltjes \cite{Stie}.
The Haussdorf moment problem is always determinate as follows
from \ref{350}  and the Stieltjes-Perron inversion formula,
see Proposition \ref{350}.

For a nice overview of the early history of the moment
problem, see Kjeldsen \cite{Kjel}.
\end{sspar}

\newpage 
\subsection{Jacobi operators}

\begin{sspar}{\label{3300}} A tridiagonal matrix of the form
$$
J =\begin{pmatrix}  b_0 & a_0 & 0   & 0   & 0 & 0 & \ldots\cr
             a_0 & b_1 & a_1 & 0   & 0 & 0 & \ldots\cr
             0  & a_1 & b_2 & a_2  & 0 & 0 & \ldots\cr
             0  & 0   & a_2 & b_3& a_3 & 0 & \ldots\cr
             \vdots&  &\ddots &\ddots &\ddots &\ddots&
   \end{pmatrix}
$$
is a Jacobi operator, or an infinite Jacobi matrix,
if $b_i\in\R$ and $a_i>0$.  If $a_i=0$ for some $i$,
the Jacobi matrix splits as the direct sum of two Jacobi
matrices, of which the first is an $(i+1)\times(i+1)$-matrix.

We consider $J$ as an operator defined on the
Hilbert space $\Hp$, see Example \ref{210}. So with respect
to the standard orthonormal basis $\{ e_k\}_{k\in\Zp}$ of $\Hp$
the Jacobi operator is defined as
\begin{equation}
J\, e_k = \begin{cases} a_k\, e_{k+1} + b_k\, e_k +
a_{k-1} \, e_{k-1}, & k\geq 1, \\
a_0\, e_1 + b_0\, e_0, & k=0. \end{cases}
\label{eq3101}
\end{equation}
Note the similarity with Theorem \ref{380}. So to
each probability measure on $\R$ with finite moments we
associate a Jacobi operator on the Hilbert space $\Hp$
from the three-term recurrence relation for the corresponding
orthonormal polynomials. However, some care is necessary,
since (\ref{eq3101}) might not define a bounded operator
on $\Hp$.

{}From (\ref{eq3101}) we extend $J$ to an operator defined on
${\mathcal D}(\Zp)$, the set of finite linear combinations of
the elements $e_k$ of the orthonormal basis of $\Hp$.
The linear subspace ${\mathcal D}(\Zp)$ is dense in $\Hp$.
{}From (\ref{eq3101}) it follows that
\begin{equation}
\langle J\,v , w\rangle = \langle v, J\, w\rangle, \qquad
\forall\, v,w\in {\mathcal D}(\Zp),
\label{eq3110}
\end{equation}
so that $J$ is a densely defined symmetric operator, see
\ref{2110}.
In particular, if $J$ is bounded on ${\mathcal D}(\Zp)$, $J$
extends to a bounded self-adjoint operator by continuity.
\end{sspar}

\begin{lem} \begin{tsspar}\label{3310}
$e_k=P_k(J)e_0$ for some polynomial $P_k$ of degree $k$
with real coefficients.
In particular, $e_0$ is a cyclic vector
for the action of $J$, i.e. the linear subspace
$\{ J^ke_0\mid k\in\Zp\}$ is dense in $\Hp$. \end{tsspar}
\end{lem}

\begin{proof} It suffices to show that $e_k=P_k(J)e_0$ for
some polynomial of degree $k$, which follows easily from
(\ref{eq3101}) and $a_i>0$ using induction on $k$.
\end{proof}


\begin{lem} \begin{tsspar}\label{3320}
If the sequences $\{a_k\}$ and
$\{b_k\}$ are bounded, say
$\sup_k |a_k| +\sup_k |b_k| \leq M <\infty$, then
$J$ extends to a bounded self-adjoint operator with
$\| J\| \leq 2M$.
On the other hand, if $J$ is bounded, then the sequences
$\{ a_k\}$ and $\{ b_k\}$ are bounded. \end{tsspar}
\end{lem}

\begin{proof} If $\{a_k\}$, $\{b_k\}$ are bounded, then,
with $v=\sum_{k=0}^\infty v_ke_k \in {\mathcal D}(\Zp)$,
$\|v\|=1$,
$$
\| Jv\|^2 = \sum_{k=0}^\infty | a_k v_{k-1} +b_k v_k
+a_{k-1}v_{k+1}|^2.
$$
Let $\sup_k |a_k|=A$, $\sup_k |b_k|=B$ with $A+B\leq M$, then
each summand can be written as
\begin{equation*}
\begin{split}
&a_k^2 |v_{k-1}|^2 +b_k^2 |v_k|^2
+a_{k-1}^2|v_{k+1}|^2 +2a_kb_k\Re(v_{k-1}\overline{v_k})
+2a_ka_{k-1}\Re(v_{k-1}\overline{v_{k+1}})
 +2b_ka_k\Re(v_{k+1}\overline{v_k}) \\ \leq&
A^2(|v_{k-1}|^2 + |v_{k+1}|^2) + B^2|v_k|^2 +
2A^2|\Re(v_{k-1}\overline{v_{k+1}})| +
2AB|\Re(v_{k-1}\overline{v_k})|
+2AB|\Re(v_{k+1}\overline{v_k})|.
\end{split}
\end{equation*}
Using the bounded shift operator $S\colon e_k\mapsto e_{k+1}$ we
have, using the Cauchy-Schwarz inequality \ref{205} and
$\| S\| =1$,
\begin{equation*}
\begin{split}
\| Jv\|^2 &\leq 2A^2+B^2 + 2A^2|\langle S^2v,v\rangle| +
4AB|\langle Sv,v\rangle| \leq 4A^2 + B^2 +4AB \\ &= (A+B)^2
+2AB + 3A^2 = 2(A+B)^2 +2A^2-B^2 \leq 4M^2.
\end{split}
\end{equation*}
By continuity $J$ extends to a bounded operator on $\Hp$.

To prove the reverse statement, we have
$|\langle J\, e_k,e_l\rangle| \leq \|J\|$ implying that
$|a_k|\leq \|J\|$ (take $l=k+1$) and $|b_k|\leq \|J\|$
(take $l=k$).
\end{proof}

\begin{sspar}{\label{3330}}
Assume $J$ is bounded, then $J$ is self-adjoint and we
can apply the spectral theorem for bounded self-adjoint
operators, see Theorem \ref{260}. Thus
$$
\langle J\, v,w\rangle = \int_\R t \, dE_{v,w}(t), \qquad
\forall \, v,w\in \Hp.
$$
In particular, we define the measure $\mu(A) = E_{e_0,e_0}(A) =
\langle E(A)e_0, e_0\rangle$. Since $E$ is a resolution of
the identity, $E(A)$ is an orthogonal projection implying that
$\mu$ is a positive Borel measure. Indeed,
$\mu(A) = \langle E(A)^2e_0,e_0\rangle = \langle E(A)e_0,
E(A)e_0\rangle \geq 0$ for any Borel set $A$.
Moreover, $E(\R)=I$, so
that $\mu$ is a probability measure. The support of $\mu$
is contained in the interval $[-\|J\|, \|J\|]$, since
$J$ is a bounded self-adjoint operator. In particular,
$\mu$ has finite moments. Hence
the spectral theorem associates to a bounded Jacobi operator
$J$ a compactly supported probability measure $\mu$.
Moreover, the spectral measure $E$ is completely determined
by $\mu$. Indeed,
\begin{equation}
\begin{split}
\langle E(A)e_k,e_l\rangle &=
\langle E(A)P_k(J)e_0,P_l(J)e_0\rangle  \\ &=
\langle P_l(J)P_k(J) E(A)e_0,e_0\rangle =
\int_A P_k(x)P_l(x)\, d\mu(x),
\end{split}
\label{eq3120}
\end{equation}
where the polynomials $P_k$
are as in Lemma \ref{3310}
using the self-adjointness of $J$
and the fact that the spectral projections commute with $J$.
\end{sspar}

\begin{thm}\begin{tsspar}\label{3340}
Let $J$ be a bounded Jacobi operator, then there exists a
unique compactly supported probability measure $\mu$
such that for any polynomial $P$ the map
$U\colon P(J)e_0\mapsto P$ extends to a unitary operator
$\Hp\to L^2(\mu)$ with $UJ=MU$, where
$M\colon L^2(\mu)\to L^2(\mu)$ is the multiplication
operator $(Mf)(x)=xf(x)$.
Moreover, let $p_k=Ue_k$, then the set $\{p_k\}_{k=0}^\infty$
is the set of orthonormal
polynomials with respect to $\mu$;
$$
\int_\R p_k(x)p_l(x)\, d\mu(x) = \de_{k,l}.
$$
\end{tsspar}
\end{thm}

\begin{proof} By Lemma \ref{3310} we see that $U$ maps
a dense subspace of $\Hp$ onto a dense subspace of $L^2(\mu)$,
since the polynomials are dense in $L^2(\mu)$ because $\mu$ is
compactly supported. Using (\ref{eq3120}) we see that
for any two polynomials $P$, $Q$ we have
\begin{equation*}
\begin{split}
\langle P(J)e_0,Q(J)e_0\rangle &=
\langle \bar Q(J)P(J)e_0,e_0\rangle
= \int_\R \bar Q(x) P(x)\, d\mu(x)
\\ &= \langle P, Q\rangle_{L^2(\mu)}
= \langle UP(J)e_0,UQ(J)e_0\rangle_{L^2(\mu)},
\end{split}
\end{equation*}
or $U$ is unitary, and it extends uniquely to a unitary operator.

To show that $U$ intertwines the Jacobi operator $J$ with the
multiplication operator we show that $UJ\, e_k=MU\, e_k$
for all $k\in\Zp$. Define $p_k=Ue_k$,
then $p_k\in L^2(\mu)$ is a polynomial
of degree $k$ and the set $\{p_k\}_{k=0}^\infty$
is the set of orthonormal
polynomials with respect to $\mu$;
$$
\int_\R p_k(x)p_l(x)\, d\mu(x) = \langle Ue_k,
Ue_l\rangle_{L^2(\mu)} = \langle e_k,
e_l\rangle = \de_{k,l}
$$
by the unitarity of $U$. It now suffices to show that
the coefficients in the three-term recurrence
relation for the orthogonal polynomials as in Theorem \ref{380}
correspond to the matrix coefficients of $J$. This is immediate
using the functional calculus of \ref{270},
(\ref{eq3120}) and the explicit expressions for
the coefficients in the three-term recurrence relation given in
the proof of Theorem \ref{380};
\begin{equation*}
\begin{split}
a_k &= \langle J\, e_k, e_{k+1}\rangle =
\int_\R x \, dE_{e_k,e_{k+1}}(x) = \int_\R xp_k(x)p_{k+1}(x)\,
d\mu(x),\\
b_k &= \langle J\, e_k, e_k\rangle =
\int_\R x \, dE_{e_k,e_k}(x) = \int_\R x\bigl( p_k(x)\bigr)^2\,
d\mu(x).
\end{split}
\end{equation*}

To show uniqueness, we observe that the moments of $\mu$
are uniquely determined by the fact that $\{ p_k\}_{k=0}^\infty$
is a set of orthonormal polynomials for $L^2(\mu)$.
Since the measure is compactly supported
its Stieltjes transform is analytic in a neighbourhood of
$\infty$. Using the Stieltjes inversion formula of
Proposition \ref{350}, we see that the  compactly supported
measure is uniquely determined by its moments.
\end{proof}

Theorem \ref{3340} is called Favard's theorem restricted
to the case of bounded coefficients in the three-term recurrence
operator. It states that any set $\{p_k\}_{k=0}^\infty$
of polynomials generated by (\ref{eq330}), (\ref{eq340}) with the
initial condition $p_0(x)=1$ with $a_k>0$, $b_k\in\R$ are
orthonormal polynomials with respect to some, see \S 3.2,
positive probability measure $\mu$.

\begin{cor}
The moment generating function for $\mu$ is in terms of
the resolvent $R(z)$ for the Jacobi operator $J$;
$$
\int_\R \frac{d\mu(x)}{x-z} = \langle R(z)e_0, e_0\rangle
=\langle (J-z)^{-1} e_0, e_0\rangle, \qquad z\in\C\backslash\R.
$$
\end{cor}

\begin{sspar}{\label{3350}}
The asymptotically free solution to $J\, f(z)=z\, f(z)$ for
$z\in\C\backslash\R$, is the
element $f(z)=\{ f_k(z)\}_{k=0}^\infty$ satisfying
$\bigl(Jf(z)\bigr)_k =z\, f_k(z)$ for $k\geq 1$ (and in general
not for $k=0$)
and $\sum_{k=0}^\infty |f_k(z)|^2<\infty$.
The asymptotically free solution encodes
the spectral measure of $J$, hence the probability measure
$\mu$ of Theorem \ref{3340}.

\begin{prop} Let $J$ be a bounded Jacobi matrix.
Take $z\in\C\backslash\R$
fixed. Then $f(z)=(J-z)^{-1}e_0$ is the asymptotically free solution
for the Jacobi operator $J$. There exists
a unique $w(z)\in\C\backslash\R$
such that $f_k(z)=w(z)p_k(z) + r_k(z)$, with $p_k$, $r_k$
the polynomials as in {\rm \ref{380}}, {\rm \ref{3100}}.
Moreover, $w$ is the Stieltjes
transform of the measure $\mu$.
\end{prop}

\begin{proof} We have already observed in Corollary \ref{3100}
that $\{ w(z)p_k(z) + r_k(z)\}_{k=0}^\infty\in\Hp$.
Next we consider
$$
f_k(z) = \langle (J-z)^{-1}e_0,e_k\rangle =
\int_\R \frac{p_k(x)}{x-z}\, d\mu(x) = w(z)p_k(z) + r_k(z)
$$
again by \ref{3100}. Hence, $\bigl( Jf(z)\bigr)_k = z\, f_k(z)$
for $k\geq 1$.

It remains to show uniqueness. If not, then there would be
two linearly independent solutions, so that we could combine
to get $\sum_{k=0}^\infty |p_k(z)|^2<\infty$ for
$z\in\C\backslash\R$, so that $J$ would have a non-real eigenvalue
$z$ contradicting the self-adjointness of $J$.
\end{proof}
\end{sspar}

\begin{sspar}{\label{3355}}
Note that we need $\lim_{K\to\infty}w(z)p_k(z)+r_k(z)=0$
in order to have $\{ f_k(z)\}_{k=0}^\infty \in\Hp$.
This implies
$$
w(z) = -\lim_{k\to\infty} \frac{r_k(z)}{p_k(z)}, \qquad
z\in\C\backslash\R,
$$
assuming that the limit in the right hand side exists. The
fraction $\frac{r_k(z)}{p_k(z)}$ has no non-real poles due
to the following lemma.

\begin{lem} The zeroes of $p_k(x)$ are real and simple.
\end{lem}

\begin{proof} Define
$$
J_N = \begin{pmatrix}  b_0 & a_0 & 0   & 0    &\ldots & 0 \cr
                a_0 & b_1 & a_1 & 0    &\ldots & 0 \cr
              \vdots &  &\ddots &\ddots& & \vdots \cr
              0 &\ldots  & & a_{N-2}   & b_{N-1} &  a_{N-1} \cr
                0 &\ldots  & & 0   & a_{N-1} &  b_N
\end{pmatrix},
$$
i.e. a tridiagonal matrix that is obtained from $J$ by
keeping only the first $(N+1)\times (N+1)$ block matrix.
Then $J_N=J_N^\ast$, and it follows that its spectrum is real.
Moreover, its spectrum is simple. Indeed, $(J_N-\la)f=0$ and
$f_0=0$ implies $f_1=0$ and hence $f_k=0$. So if the multiplicity
of the eigenspace is more than one, we could construct a
non-zero eigenvector with $f_0=0$, a contradiction.

On the other hand we have, from Theorem \ref{3340},
$$
(J_N -z) \begin{pmatrix} p_0(z) \cr p_1(z) \cr \vdots \cr p_N(z)
\end{pmatrix} = -a_N \begin{pmatrix} 0\cr\vdots\cr 0\cr p_{N+1}(z)
\end{pmatrix},
$$
so that the eigenvalues of $J_N$ are the zeroes
of $p_{N+1}$.
\end{proof}
\end{sspar}

\begin{sspar}{\label{3360}}
The norm of the asymptotically free solution for fixed
$z\in\C\backslash\R$ can be expressed in terms of the
Stieltjes transform of $\mu$;
\begin{equation}
\begin{split}
&\sum_{k=0}^\infty |w(z)p_k(z)+r_k(z)|^2 = \langle (J-z)^{-1}e_0,
(J-z)^{-1}e_0\rangle  \\ &= \frac{1}{z-\bar z}\bigl(
\langle (J-z)^{-1}e_0,e_0\rangle -
\langle (J-\bar z)^{-1}e_0,e_0\rangle \bigr) =
\frac{w(z)-\overline{w(z)}}{z-\bar z}.
\end{split}
\label{eq3165}
\end{equation}
For $z$ fixed and $w=w(z)$ a complex parameter (\ref{eq3165})
gives rise to an equation in the complex $w$-plane, which is
in general a circle or a point. The radius of this circle
is $\bigl( |z-\bar z|\sum_{k=0}^\infty |p_k(z)|^2\bigr)^{-1}$.
Proposition \ref{3350} shows that it is
a point for a bounded Jacobi operator, see \cite{Akhi},
\cite{Simo} for a discussion of limit points and limit circles.
\end{sspar}

\begin{sspar}{\label{3370}}
We now introduce the Green kernel for $z\in\C\backslash\R$,
$$
G_{k,l}(z) = \begin{cases} f_l(z)p_k(z), & k\leq l, \\
                    f_k(z)p_l(z), & k>l. \end{cases}
$$
Hence $\{ G_{k,l}(z)\}_{k=0}^\infty,
\{ G_{k,l}(z)\}_{l=0}^\infty \in \Hp$ by
Proposition \ref{3350}, and the map
$$
\Hp \ni v \mapsto (Gz)v, \qquad
\bigl( G(z)v\bigr)_k = \sum_{l=0}^\infty v_l G_{k,l}(z) =
\langle v, \overline{G_{k,\cdot}(z)}\rangle
$$ is
a well-defined map. A priori it is not clear that $G(z)$ is a
bounded map, but it is densely defined, e.g. on ${\mathcal D}(\Zp)$.
The next proposition states that $G$ is bounded.

\begin{prop}
The resolvent of $J$ is given by
$(J-z)^{-1}=G(z)$ for $z\in\C\backslash\R$.
\end{prop}

\begin{proof} Since $z\in\C\backslash\R\subset \rho(J)$ we
know that $(J-z)^{-1}$ is a bounded operator, so it suffices
to check $(J-z)G(z)={\mathbf 1}_{\Hp}$
on a dense subspace ${\mathcal D}(\Zp)$ of $\Hp$. Now
\begin{equation*}
\begin{split}
\bigl( (J-z)G(z)v\bigr)_k &= \sum_{l=0}^\infty
v_l\bigl( a_kG_{k+1,l}(z) + (b_k-z)G_{k,l}(z) +
a_{k-1}G_{k-1,l}(z)\bigr) \\
& = \sum_{l=0}^{k-1} v_lp_l(z)
\bigl(a_kf_{k+1}(z) + (b_k-z)f_k(z) +a_{k-1}f_{k-1}(z)\bigr)  \\
& \qquad  +\sum_{l=k+1}^\infty v_lf_l(z)
\bigl(a_kp_{k+1}(z) + (b_k-z)p_k(z) +a_{k-1}p_{k-1}(z)\bigr)  \\
&  \qquad  +v_k\bigl( a_k f_{k+1}(z)p_k(z) + (b_k-z) f_k(z)p_k(z)
+ a_{k-1} f_k(z) p_{k-1}(z) \bigr)\\
& = v_k\bigl( a_k f_{k+1}(z)p_k(z) + (b_k-z) f_k(z)p_k(z)
+ a_{k-1} f_k(z) p_{k-1}(z) \bigr) \\
& = v_k a_k \bigl( f_{k+1}(z)p_k(z) - f_k(z)p_{k+1}(z)\bigr)
\end{split}
\end{equation*}
using that $f_k(z)$ and $p_k(z)$ are solutions to $Jf=zf$ for
$k\geq 1$. The term at the right hand side involves the Wronskian,
see \ref{3150}. Since $[p,p]=0$, and $[p,r]=-1$ and
$f_k(z)=w(z)p_k(z)+r_k(z)$, it follows
that $[f,p]=
a_k \bigl( f_{k+1}(z)p_k(z) - f_k(z)p_{k+1}(z)\bigr)= 1$.
\end{proof}
\end{sspar}

\begin{sspar}{\label{3380}}
We can now find the spectral measure of $J$ from the resolvent,
see the Stieltjes-Perron inversion formula \ref{280};
$$
E_{u,v}\bigl( (a,b)\bigr) = \lim_{\de\downarrow 0}
\lim_{\ep\downarrow 0} \frac{1}{2\pi i}
\int_{a+\de}^{b-\de}
\langle G(x+i\ep)u,v\rangle - \langle G(x-i\ep)u,v\rangle \, dx.
$$
First observe that
\begin{equation*}
\begin{split}
\langle G(z)u,v\rangle = \sum_{k,l=0}^\infty
G_{k,l}(z) u_l \bar v_k &=
\sum_{k\leq l} f_l(z)p_k(z) u_l \bar v_k +
\sum_{k> l}f_k(z)p_l(z) u_l \bar v_k \\ &=
\sum_{k\leq l} f_l(z)p_k(z) (u_l \bar v_k + u_k\bar v_l)(1-\hf
\de_{k,l}),
\end{split}
\end{equation*}
by splitting the sum and renaming summation variables. The
factor $(1-\hf\de_{k,l})$ is introduced in order to avoid
doubling in the case $k=l$. Since $f_l(z)=w(z)p_l(z) + r_l(z)$,
with $p_l$ and $r_l$ polynomials,
the only term contributing to
$\langle G(x+i\ep)u,v\rangle - \langle G(x-i\ep)u,v\rangle$
as $\ep\downarrow 0$ comes from $w(z)$, cf. \ref{350}.
Hence,
\begin{equation*}
\begin{split}
&E_{u,v}\bigl( (a,b)\bigr) = \\ &\lim_{\de\downarrow 0}
\lim_{\ep\downarrow 0} \frac{1}{2\pi i}
\int_{a+\de}^{b-\de} \Bigl( w(x+i\ep)-w(x-i\ep)\Bigr)
\sum_{k\leq l} p_l(x)p_k(x) (u_l \bar v_k + u_k\bar v_l)(1-\hf
\de_{k,l})\, dx
\end{split}
\end{equation*}
and using Proposition \ref{350} and symmetrising the sum gives
\begin{equation}
E_{u,v}\bigl( (a,b)\bigr) = \int_{(a,b)} \bigl( Uu\bigr)(x)
\overline{\bigl( Uv\bigr)(x)}\, d\mu(x),
\label{eq3185}
\end{equation}
where $U\colon \Hp\to L^2(\mu)$ is the unitary operator
$U$ of Theorem \ref{3340}. Note that (\ref{eq3185})
proves that $U$ is a unitary operator by letting $a\to -\infty$ and
$b\to \infty$, so that
$$
\langle u,v\rangle = \int_\R \bigl( Uu\bigr)(x)
\overline{\bigl( Uv\bigr)(x)}\, d\mu(x) = \langle Uu,
Uv\rangle_{L^2(\mu)}.
$$
We can think of $U$ as the Fourier
transform;
$$
u_k = \int_\R \bigl( Uu\bigr)(x) p_k(x)\, d\mu(x),
$$
so that $u\in\Hp$ is expanded in terms of (generalised)
eigenvectors $\{ p_k(x)\}_{k=0}^\infty$ of the Jacobi
operator $J$.
\end{sspar}


\begin{thm} \begin{tsspar}\label{3390} There is one-to-one
correspondence between bounded Jacobi operators and
probability measures on $\R$ with compact support.
\end{tsspar}\end{thm}

\begin{proof} Since a probability measure $\mu$
of compact support has
finite moments, we can build the corresponding orthonormal
polynomials and write down the three-term recurrence relation
by Theorem \ref{380}. This gives a map $\mu\mapsto J$.
Using the moment generating function, analytic at $\infty$,
and the Stieltjes inversion formula, cf. proof of Theorem
\ref{3340}, we see that the map is injective. In
Theorem \ref{3340} the inverse map has been constructed.
\end{proof}

\subsection{Unbounded Jacobi operators}

\begin{sspar}{\label{3400}}
We now no longer assume that $J$ is bounded. This occurs for example
in the following case.

\begin{lem} Let $\mu$ be a probablity measure with finite
moments. Consider the three-term recurrence relation of
Theorem {\rm \ref{380}} and the corresponding densely defined
Jacobi operator $J$. If $\text{\rm supp}(\mu)$ is unbounded, then
$J$ is unbounded.
\end{lem}

\begin{proof} Suppose not, so we assume that $J$ is bounded. Using
the Cauchy-Schwarz inequality and (\ref{eq3120}) we obtain
\begin{equation*}
\begin{split}
\|J\|^{2k} &\geq |\langle J^{2k}e_0,e_0\rangle | =
\big\vert \int_\R x^{2k}\, d\mu(x)\big\vert
\\ &\geq \int_{|x|\geq \|J\|+1} x^{2k}\, d\mu(x)
\geq \bigl( |J\|+1\bigr)^{2k} \int_{|x|\geq \|J\|+1} \, d\mu(x).
\end{split}
\end{equation*}
This implies
$$
\Bigl( \frac{\|J\|}{\|J\| +1}\Bigr)^{2k} \geq
\int_{|x|\geq \|J\|+1} \, d\mu(x)
$$
and by the assumption on $\mu$ the right hand side is
strictly positive. The left hand side tends to $0$ as $k\to\infty$,
so that we obtain the required contradiction.
\end{proof}
\end{sspar}

\begin{sspar}{\label{3410}}
We use the notions as recalled in \S 2.3.
Recall from (\ref{eq3110}) that $J$ is a densely defined symmetric
operator.
Let us extend $J$ to an arbitrary vector
$v=\sum_{k=0}^\infty v_ke_k\in\Hp$ by formally
putting
\begin{equation}
J^\ast \, v = (a_0v_1+b_0v_0)\, e_0 + \sum_{k=1}^\infty
\bigl( a_kv_{k+1}+b_kv_k+a_{k-1}v_{k-1}\bigr)\, e_k,
\label{eq3412}
\end{equation}
which, in general, is not an element of $\Hp$.

\begin{prop}
The adjoint of $(J, {\mathcal D}(\Zp))$ is $(J^\ast, {\mathcal D}^\ast)$,
where
$$
{\mathcal D}^\ast = \{ v\in \Hp\mid J^\ast v\in \Hp\}
$$
and $J^\ast$ is given by \textup{(\ref{eq3412})}.
\end{prop}

Proposition \ref{3410}
says that the adjoint of $J$ is its natural extension to
its maximal domain. In case $J$ is essentially self-adjoint
we have that $(J^\ast, {\mathcal D}^\ast)$ is self-adjoint and that
$(J^\ast, {\mathcal D}^\ast)$ is the closure of $J$.

\begin{proof} To determine the domain of $J^\ast$ we
have to consider for which $v\in\Hp$
the map $w\mapsto \langle Jw, v\rangle$, ${\mathcal D}(\Zp)\to \C$
is continuous, see \ref{290}. Now, with the convention
$e_{-1}=0$,
\begin{equation*}
\begin{split}
|\langle Jw, v\rangle| &= \big\vert \sum_k w_k
\langle a_k e_{k+1} +b_ke_k +a_{k-1}e_{k-1},v\rangle\big\vert \\
&=\big\vert \sum_k w_k (a_k \bar v_{k+1} +b_k\bar v_k +a_{k-1}
\bar v_{k-1}) \big\vert \\ &\leq \| w\|
\Bigl( \sum_k |a_k v_{k+1} +b_k v_k +a_{k-1}v_{k-1}|^2\Bigr)^\hf =
\| w\|  \| J^\ast v\|,
\end{split}
\end{equation*}
where the sums are finite since $w\in {\mathcal D}(\Zp)$ and where
we have used the Cauchy-Schwarz inequality. This proves
$w\mapsto \langle Jw, v\rangle$ is continuous for $v\in {\mathcal
D}^\ast$. Hence, we have proved
${\mathcal D}^\ast \subseteq {\mathcal D}(J^\ast)$.

For the other inclusion, we observe that for $v\in{\mathcal D}(J^\ast)$
we have
$$
\big\vert \sum_k w_k (a_k \bar v_{k+1} +b_k\bar v_k +a_{k-1}
\bar v_{k-1}) \big\vert \leq C \| w\|
$$
for some constant $C$ independent of $w\in{\mathcal D}(\Zp)$.
Specialise to
$w_k = (Jv)_k$ for $0\leq k\leq N$ and $w_k=0$ for $k>N$, to
find $\| w\|^2\leq C\| w\|$, or
$$
\Bigl( \sum_{k=0}^N
|a_k v_{k+1} +b_k v_k +a_{k-1}v_{k-1}|^2\Bigr)^\hf \leq C
$$
and since $C$ is independent of $w$, $C$ is independent of
$N$. Letting $N\to\infty$ we find $\|J^\ast v\|\leq C$, or
$v\in {\mathcal D}^\ast$.
\end{proof}

Since the coefficients are real, $J^\ast$ commutes with
complex conjugation. This implies that the deficiency indices
are equal. Since a solution of $J^\ast\, f=z\, f$ is completely
determined by $f_0=\langle f, e_0\rangle$, it follows that the
deficiency spaces are either zero-dimensional
or one-dimensional. In
the first case we see that $J$ with domain ${\mathcal D}(\Zp)$
is essentially self-adjoint, and in the second case there exists
a one-parameter family of self-adjoint extensions of
$(J, {\mathcal D}(\Zp))$.
Since the only possible element in $N_z$ is $\{
p_k(z)\}_{k=0}^\infty$ we obtain the following
corollary to Proposition \ref{2120}.

\begin{cor} $(J, {\mathcal D}(\Zp))$
is essentially self-adjoint if and only if
$\sum_{k=0}^\infty |p_k(z)|^2 =\infty$ for all $z\in\C\backslash\R$
if and only if
$\sum_{k=0}^\infty |p_k(z)|^2 =\infty$ for some $z\in\C\backslash\R$.
\end{cor}
\end{sspar}

\begin{sspar}{\label{3420}}
In case $J$ has deficiency indices $(1,1)$ every self-adjoint
extension, say $(J_1,{\mathcal D}_1)$,
satisfies $(J,{\mathcal D}(\Zp))\subsetneq
(J_1, {\mathcal D}_1)\subsetneq (J^\ast, {\mathcal D}^\ast)$, and by
Proposition \ref{3410} $(J_1,{\mathcal D}_1)$ is the restriction of
$J^\ast$ to a domain ${\mathcal D}_1$. Since the polynomials
$p_k$ and $r_k$ are generated
from the restriction to ${\mathcal D}(\Zp)$, these polynomials
are independent of the choice of the self-adjoint extension.
To describe the domains of the self-adjoint extensions
we rewrite the sesquilinear form $B$ of \ref{2120}
in terms of the Wronskian; using the convention
$u_{-1}=0=v_{-1}$
\begin{equation*}
\begin{split}
&\sum_{k=0}^N (J^\ast u)_k \bar v_k - u_k \overline{(J^\ast v)_k}
\\ &=
\sum_{k=0}^N (a_ku_{k+1}+b_ku_k+a_{k-1}u_{k-1})\bar v_k
-u_k(a_k\bar v_{k+1}+b_k\bar v_k+a_{k-1}\bar v_{k-1}) \\ &=
[u,\bar v]_0 + \sum_{k=1}^N \Bigl(
[u,\bar v]_k - [u,\bar v]_{k-1}\Bigr) = [u,\bar v]_N,
\end{split}
\end{equation*}
where $[u,v]_k =a_k (u_{k+1}v_k-v_{k+1}u_k)$ is the
Wronskian, cf. (\ref{eq351}). So for $u,v\in {\mathcal D}^\ast$ we have
$B(u,v)=\lim_{N\to\infty} [u,v]_N$.
Note that the limit exists, since $u,v\in {\mathcal D}^\ast$.

\begin{lem} Assume that
$J$ has deficiency indices $(1,1)$, then
the self-adjoint extensions are in one-to-one
correspondence with $(J^\ast, {\mathcal D}_\te)$,
$\te\in [0,2\pi)$,  where
$$
{\mathcal D}_\te = \{ v\in {\mathcal D}^\ast \mid
\lim_{N\to\infty} [v, e^{i\te}\psi_i + e^{-i\te}\psi_{-i}]_N=0\}
$$
where $J^\ast \psi_{\pm i} = \pm i\, \psi_{\pm i}$,
$\overline{(\psi_i)_k} = (\psi_{-i})_k$ and
$\|\psi_{\pm i}\|=1$.
\end{lem}

So $(\psi_{\pm i})_k =C\, p_k(\pm i)$ with
$C^{-2}=\sum_{k=0}^\infty |p_k(i)|^2<\infty$ by Corollary
\ref{3410}.

\begin{proof} Since $n_+=n_-=1$, Proposition \ref{2120} gives
that the domain of a self-adjoint extension is of the
form $u+c(\psi_i+e^{-2i\te}\psi_{-i})$ with $u\in {\mathcal
D}(J^{\ast\ast})$, $c\in\C$. Using $B(\psi_i,\psi_i)=2i$,
$B(\psi_{-i},\psi_{-i})=-2i$ and $B(\psi_i,\psi_{-i})=0$
shows that
\begin{equation*}
\begin{split}
&B(u+c(\psi_i+e^{-2i\te}\psi_{-i}), e^{i\te}\psi_i +
e^{-i\te}\psi_{-i}) = B(u, e^{i\te}\psi_i +
e^{-i\te}\psi_{-i}) \\ & =
\langle J^\ast u, e^{i\te}\psi_i +
e^{-i\te}\psi_{-i}\rangle  - \langle u, e^{i\te}J^\ast \psi_i +
e^{-i\te}J^\ast \psi_{-i}\rangle.
\end{split}
\end{equation*}
Now observe that for $v\in N_{\pm i}$ we have
$B(u,v)= {\pm i} \langle u,
v\rangle_{J^\ast}=0$, $u\in {\mathcal D}(J^{\ast\ast})$, using
the graph norm as in
Proposition \ref{2120}.
\end{proof}

The condition $\overline{(\psi_i)_k} = (\psi_{-i})_k$ in
Lemma \ref{3420} is meaningful, since $J^\ast$ commutes
with complex conjuagtion. It is needed to ensure the
one-to-one correspondence.
\end{sspar}

\begin{sspar}{\label{3430}} The results of \S 3.3 go through in the
setting of unbounded operators up to minor changes.
Lemma \ref{3310} remains valid, as it has nothing to do
with the unboundedness of $J$. Lemma \ref{3320} implies
that for an unbounded Jacobi matrix we have to deal with
unbounded sequences $\{a_k\}$, $\{b_k\}$, or at least one
of them is unbounded.
\ref{3330} goes through, using Theorem  \ref{2130},
after choosing a self-adjoint extension of $J$. This operator
is uniquely determined
if $(n_+,n_-)=(0,0)$ and is labeled by one real parameter
if $(n_+,n_-)=(1,1)$. For this self-adjoint extension
\ref{3330} and
Theorem \ref{3340} remain valid, except that $\mu$ is no longer
compactly supported. In this case we observe that $\mu$ has
finite moments, since $\int_\R x^k\, d\mu(x) =
\langle J^ke_0,e_0\rangle<\infty$. From the unicity
statement in Theorem \ref{2130} we deduce
that the polynomials are dense in $L^2(\mu)$, which is needed
to obtain the unitarity of $U$ in Theorem \ref{3340}.
The results of \S 3.3
until Theorem \ref{3390} remain valid after replacing $J$ by a
self-adjoint extension. Observe that $f$ as in Proposition
\ref{3350} is contained in the domain of a self-adjoint
extension, so that the calculation in Proposition \ref{3370}
remains valid.
In the unbounded case it is no
longer true that (\ref{eq3165}) in the complex $w$-plane
describes a point.
Theorem \ref{3440}
and Corollary \ref{3410} show that in the
general case (\ref{eq3165}) describes a circle (with positive radius)
if and only if the Jacobi operator has deficiency indices $(1,1)$
if and only if the corresponding (Hamburger) moment problem is
indeterminate. Finally,
Theorem \ref{3390} cannot be extended
to unbounded self-adjoint Jacobi operators.
\end{sspar}

\begin{sspar}{\label{3440}} There is a nice direct link between the
moment problem in \S 3.2 and the Jacobi operators. It gives
an answer to question 2 of \ref{3250}.

\begin{thm} The (Hamburger) moment problem is
determinate if and only if the corresponding Jacobi operator
is essentially self-adjoint.
\end{thm}

\begin{proof} Assume first that the corresponding Jacobi operator
is essentially self-adjoint. Then $(J-z)\bigl({\mathcal D}(\Zp)\bigr)$
is dense in $\Hp$ by Proposition \ref{2120},
$z\in\C\backslash\R$. The density and
Lemma \ref{3310} say that we can find polynomials
$R_k(x;z)$ such that $\|(J-z)R_k(J;z)e_0-e_0\|\to 0$ as $k\to\infty$.
Let $\mu$ be any solution to the moment problem, so that
$J$ is realised as the multiplication operator on the closure of
the span of the polynomials in $L^2(\mu)$. Then
$$
\lim_{k\to\infty} \int_\R |(x-z)R_k(x;z)-1|^2\, d\mu(x) =0.
$$
Since $z\in\C\backslash\R$ we see that $x\mapsto (x-z)^{-1}$
is bounded for $x\in\R$, hence
$$
\lim_{k\to\infty} \int_\R |R_k(x;z)-\frac{1}{x-z}|^2\, d\mu(x) =0.
$$
Since $\mu$ is a probability measure
$L^2(\mu) \hookrightarrow L^1(\mu)$, so
the Stieltjes transform of $\mu$ is given by
$$
\int_\R \frac{d\mu(x)}{x-z} = \lim_{k\to\infty} \int_\R R_k(x;z)\,
d\mu(x),
$$
but the right hand side only involves integration of polynomials,
hence only the moments, so it is independent of the choice
of the solution of the moment problem. By the
Stieltjes inversion of Proposition \ref{350} this determines
$\mu$ uniquely.

For the converse statement, we assume that $J$ has deficiency
indices $(1,1)$. Let $J_1$ and $J_2$ be two different
self-adjoint extensions of $(J,{\mathcal D}(\Zp))$. We have to
show that they give rise to different solutions of the
moment problem. By the Stieltjes inversion of Proposition
\ref{350} it suffices to show that the Stieltjes transforms
$\langle (J_1-z)^{-1} e_0,e_0\rangle$ and
$\langle (J_2-z)^{-1} e_0,e_0\rangle$ are different.

Let us first observe that $e_0\notin
\text{\rm Ran}(J^{\ast\ast}-z)$. Indeed, suppose not, then
there exists $u\in{\mathcal D}(J^{\ast\ast})$ such that
$(J^{\ast\ast}-z)u=e_0$, and taking the inner product with
$0\not= v\in N_{\bar z}$
gives $\langle e_0, v\rangle = \langle (J^{\ast\ast}-z)u,v\rangle =
\langle u, (J^\ast - \bar z)v\rangle = 0$,
since $u\in {\mathcal D}(J^{\ast\ast})$
and $v\in {\mathcal D}(J^\ast)$ and $v\in N_{\bar z}$. So
$v$ is an eigenvector of $J^\ast$ with $v_0=0$, so that
$v=0$ by \ref{3410}. This contradicts the fact that $J$ has
deficiency indices $(1,1)$.

{}From the previous paragraph we conclude that $(J_i-z)^{-1}e_0\in
{\mathcal D}(J^\ast)\backslash {\mathcal D}(J^{\ast\ast})$, $i=1,2$, by
Proposition \ref{2120}. Now the
$\dim \bigl( {\mathcal D}(J_i)\backslash {\mathcal D}(J^{\ast\ast})\bigr)
=1$, so that $(J_1-z)^{-1}e_0 = (J_1-z)^{-1}e_0$ implies
${\mathcal D}(J_1)={\mathcal D}(J_2)$ and hence $J_1=J_2$ by
Proposition \ref{2120}. So the vectors $(J_i-z)^{-1}e_0$
are different. To finish the proof  we need to show that the
zero-components are also different. Let
$u=(J_1-z)^{-1}e_0-(J_2-z)^{-1}e_0\not=0$, then
$$
(J^\ast-z)u=(J^\ast-z)(J_1-z)^{-1}e_0-(J^\ast-z)(J_2-z)^{-1}e_0
=e_0-e_0=0,
$$
or $u\in N_z$. By \ref{3410} $\langle u,e_0\rangle =0$ implies
$u=0$, so $\langle u,e_0\rangle \not=0$. Hence, $J_1$ and $J_2$
give rise to two different Stieltjes transforms, hence to
two different solutions of the moment problem.
\end{proof}
\end{sspar}

\begin{sspar}{\label{3450}} From Lemma \ref{3320} we see that boundedness
of the coefficients in the three-term recurrence relation
of Theorem \ref{380} for the orthonormal polynomials implies
that $J$ extends to a bounded self-adjoint operator. For unbounded
coefficients there are several conditions on the sequences
$\{a_k\}$, $\{b_k\}$ that ensure that the Jacobi operator
$J$ is essentially
self-adjoint, and hence, by Theorem \ref{3440}, the
corresponding moment problem is determinate.
We give three examples.

\begin{prop} {\rm (i)} If $\sum_{k=0}^\infty
\frac{1}{a_k}=\infty$, then $J$ is essentially self-adjoint.
\par\noindent
{\rm (ii)} If $a_k+b_k+a_{k-1}\leq M<\infty$ for $k\in\N$, or
if $a_k-b_k+a_{k-1}\leq M<\infty$ for $k\in\N$, then
$J$ is essentially self-adjoint. \par\noindent
{\rm (iii)} If $\sum_{k=0}^\infty m_{2k}^{-1/2k}=\infty$,
where $m_j$ is the $j$-th moment $m_j=\int_\R x^j\, d\mu(x)$,
then $J$ is essentially self-adjoint.
\end{prop}

\begin{proof} To prove (i), we use the Christoffel-Darboux
formula (\ref{eq360}) with $x=z$, $y=\bar z$, and using that
the polynomials have real coefficients, we find for
$z\in\C\backslash\R$,
$$
\sum_{k=0}^N |p_k(z)|^2= \frac{a_N}{z-\bar z} \bigl( p_{N+1}(z)
\bar p_N(z) - p_N(z)\bar p_{N+1}(z)\bigr).
$$
Since the left hand side is at least $1=p_0(z)$, we obtain
$$
1\leq \frac{a_N}{|z-\bar z|} \bigl( |p_{N+1}(z)|
| p_N(z)| + |p_N(z)| |p_{N+1}(z)|\bigr)
= \frac{2a_N}{|z-\bar z|} |p_N(z)| |p_{N+1}(z)|.
$$
For non-real $z$ we have $C=2/|z-\bar z|>0$ and so
\begin{equation}
\sum_{k=0}^\infty \frac{1}{a_k} \leq C
\sum_{k=0}^\infty |p_k(z)| |p_{k+1}(z)| \leq
C \sum_{k=0}^\infty |p_k(z)|^2
\label{eq3451}
\end{equation}
by the Cauchy-Schwarz inequality \ref{205} in $\Hp$.
So, if $\sum_{k=0}^\infty
\frac{1}{a_k}$ diverges, we see from Corollary \ref{3410}
and (\ref{eq3451}) that $J$ is essentially self-adjoint.

The cases in (ii) are equivalent, by switching
to the orthonormal polynomials $(-1)^np_n(-x)$ for
the measure $\tilde \mu(A)=\mu(-A)$ for any Borel set
$A\subset \R$.  We start with
an expression for $p_n(x)$ in terms of the
lower degree orthonormal polynomials;
\begin{equation}
p_n(x) = 1 + \sum_{k=0}^{n-1} \frac{1}{a_k}
\sum_{j=0}^k (x-b_j -a_j-a_{j-1})p_j(x),
\label{eq3452}
\end{equation}
with the convention $a_{-1}=0$ and $p_{-1}(x)=0$.
In order to prove (\ref{eq3452}) we start with the telescoping
series
\begin{equation*}
\begin{split}
a_{n-1}(p_n(x)-p_{n-1}(x)) =& \sum_{j=0}^{n-1}
\bigl( a_j(p_{j+1}(x)-p_j(x))- a_{j-1}(p_j(x)-p_{j-1}(x))\bigr)
\\ =&
\sum_{j=0}^{n-1} (x-b_j-a_j-a_{j-1})\, p_j(x)
\end{split}
\end{equation*}
by (\ref{eq330}), (\ref{eq340}). This gives
$$
p_n(x) = p_{n-1}(x) +\frac{1}{a_{n-1}}
\sum_{j=0}^{n-1} (x-b_j-a_j-a_{j-1})\, p_j(x)
$$
and iterating this expression gives (\ref{eq3452}).

Now take $x>M$ in (\ref{eq3452}) to obtain inductively
$p_0(x)=1$, $p_k(x)\geq 1$ for $k\in \N$, so that
$\sum_{k=0}^\infty |p_k(x)|^2$ is divergent. Since
the polynomials $p_k$ have real zeroes, see Lemma \ref{3355},
we have $p_k(x)=C_k\prod_{i=1}^k (x-x_i)$ with $x_i\in\R$  which
implies that $|p_k(x+iy)|\geq |p_k(x)|$.
Hence, $\sum_{k=0}^\infty |p_k(x+iy)|^2=\infty$ and
Corollary \ref{3410} shows $J$ is essentially self-adjoint.

Finally for (iii) we note that
$$
1=\| p_k\|_{L^2(\mu)}\leq \| \mathrm{lc}(p_k)\, x^{k}\|
= \mathrm{lc}(p_k) \sqrt{m_{2k}} =
\frac{\sqrt{m_{2k}}}{a_0a_1a_2\ldots a_{k-1}}
$$
by the triangle inequality, where $\mathrm{lc}(p_k)$
denotes the leading coefficient of the polynomial
$p_k$. Hence
\begin{equation*}
m_{2k}^{-1/2k}\leq \bigl(
\frac{1}{a_0}\frac{1}{a_1}\frac{1}{a_2}\ldots
\frac{1}{a_{k-1}}\bigr)^{1/k}
\Longrightarrow
\sum_{k=0}^\infty m_{2k}^{-1/2k} \leq
\sum_{k=0}^\infty \Bigl(\prod_{i=0}^{k-1} \frac{1}{a_i}\Bigr)^{1/k}.
\end{equation*}
Using the geometric-arithmetic mean inequality,
$$
\prod_{i=1}^n a_i^{p_i} \leq \sum_{i=1}^n p_i a_i, \qquad
\sum_{i=1}^n p_i=1,
$$
we find
$$
\Bigl(\prod_{i=0}^{k-1} \frac{1}{a_i} \Bigr)^{1/k}=
(k!)^{-1/k}\Bigl(\prod_{i=0}^{k-1} \frac{i+1}{a_i}\Bigr)^{1/k}
\leq (k!)^{-1/k} \frac{1}{k} \sum_{i=0}^{k-1} \frac{i+1}{a_i}.
$$
Next use $k^k\leq e^k k!$, which can be proved by induction on $k$
and $(1+\frac{1}{k})^k\leq e$, to see that we can estimate
$$
\sum_{k=0}^\infty m_{2k}^{-1/2k} \leq \sum_{k=0}^\infty
\frac{e}{k^2} \sum_{i=0}^{k-1} \frac{i+1}{a_i} =
e \sum_{i=0}^\infty \frac{1}{a_i} \sum_{k=i}^\infty \frac{i+1}{k^2}
$$
and the inner sum over $k$ can be estimated independent of $i$.
Now the result follows from (i).
\end{proof}
\end{sspar}

\begin{sspar}{\label{3455}} A famous example  of an indeterminate
moment problem has been given by Stiel\-tjes in his
posthumously published memoir \cite{Stie}. Consider
for $\ga>0$ ---Stieltjes considered the case $\ga=1$---
$$
\int_0^\infty x^n \, e^{-\ga^2\ln^2 x}\sin (2\pi\ga^2\ln x)\, dx =
e^{\frac{(n+1)^2}{4\ga^2}}
\int_0^\infty e^{-(\ga\ln x- \frac{n+1}{2\ga})^2}
\sin (2\pi\ga^2\ln x)\, \frac{dx}{x}
$$
and put $y=\ga\ln x- \frac{n+1}{2\ga}$,
$x^{-1}dx = \ga^{-1}dy$ to see that the integral
equals
$$
\ga^{-1} e^{\frac{(n+1)^2}{4\ga^2}}  \int_{-\infty}^\infty
e^{-y^2} \sin (2\pi\ga y + \pi(n+1))\, dy = 0,
$$
since the integrand is odd. In the same way we can calculate
the same integral without the sine-term and we
find that for any $-1\leq r\leq 1$ we have
$$
\int_0^\infty x^n \, e^{-\ga^2\ln^2 x}\bigl(
1+r\sin (2\pi\ga^2\ln x)\bigr)\, dx =
\frac{\sqrt{\pi}}{\ga} \exp\bigl(\frac{(n+1)^2}{4\ga^2}\bigr),
$$
so that the probablity measure
$$
d\mu(x) = \frac{\ga}{\sqrt{\pi}} e^{\frac{-1}{4\ga^2}}\,
e^{-\ga^2\ln^2 x}\, dx, \qquad x>0,
$$
corresponds to an indeterminate moment problem.
The corresponding orthogonal polynomials are known
as the Stieltjes-Wigert polynomials. Put
$q=\exp(-1/2\ga^2)$ or $\ga^{-2}=-2\ln q$, then the
orthonormal Stieltjes-Wigert polynomials satisfy
(\ref{eq330}), (\ref{eq340}) with
$a_k=q^{-2k-\frac{3}{2}}\sqrt{1-q^{k+1}}$ and
$b_k=q^{-2k}(1+q^{-1}-q^k)$, see e.g. \cite[\S VI.2]{Chih}.
Note that $0<q<1$, and that $a_k$ and $b_k$ are exponentially
increasing.
\end{sspar}

\begin{sspar}{\label{3457}} Another example is to consider the
following measure on $\R$
$$
d\mu(x) = C_{\al,\ga} \, e^{-\ga|x|^\al}\, dx, \qquad
\al,\ga>0.
$$
The moments and the explicit value
for $C_{\al,\ga}$ can be calculated using the $\Ga$-function;
\begin{equation}
\int_0^\infty x^{c-1}e^{-bx}\, dx = b^{-c}\Ga(c),
\qquad c>0, \ \Re b>0.
\label{eq3458}
\end{equation}
For $0<\al<1$ the (Hamburger) moment problem is indeterminate
as we can see from the following observations. Note that
all odd moments vanish. Put $x=y^\al$ and $c=(2n+1)/\al$ to
find from (\ref{eq3458}), after doubling the interval,
$$
\frac{\al}{2} \int_{-\infty}^\infty y^{2n} e^{-b|y|^\al}\, dy =
b^{-\frac{2n+1}{\al}} \Ga\bigl( \frac{2n+1}{\al}\bigr).
$$
By taking the real parts we obtain
$$
\frac{\al}{2} \int_{-\infty}^\infty y^{2n} e^{-\Re b |y|^\al}
\cos (-\Im b |y|^\al) \, dy = \Ga\bigl( \frac{2n+1}{\al}\bigr)
|b|^{-\frac{2n+1}{\al}} \cos (- \arg(b) \frac{2n+1}{\al}),
$$
so that the right hand side is zero for all $n\in\Zp$ if
$\arg (b) = \hf \al\pi$. Since we need $\Re b>0$ this
is possible if $|\al|<1$.  So for $0<\al<1$ the (Hamburger)
moment problem is indeterminate, since we obtain more solutions
in the same way as for the Stieltjes-Wigert polynomials
in \ref{3455}.

The (Hamburger) moment problem is determinate for $\al\geq 1$,
which follows using Proposition \ref{3450}(iii).
See also Deift \cite[p.~34]{Deif} for another
proof, in which the essentially
self-adjointness of the corresponding
Jacobi operator is established.

Restricting the measure $d\mu(x)$ to $[0,\infty)$ gives
a Stieltjes moment problem, see \S 3.2, and for the Stieltjes
moment problem this is an indeterminate moment problem
if and only if $0<\al<\hf$, see \cite[Ch.~1, \S 8]{ShohT}.
In particular, for $\hf\leq \al<1$, the Hamburger moment
problem is indeterminate, whereas the Stieltjes moment
problem is determinate. In terms of Jacobi operators, this
means that the corresponding Jacobi operator $J$ has
deficiency indices $(1,1)$, but that $J$ has a unique
positive self-adjoint extension, see Simon \cite{Simo}.
\end{sspar}

\begin{sspar}{\label{3460}} A great number of results exist for
the description of the solutions to
an indeterminate moment problem. A nice way to describe all
solutions is due to Nevannlina; see e.g. \cite{Akhi},
\cite{Simo}.
There exist four entire functions $A$, $B$, $C$ and $D$,
that can be described using the polynomials $p_k$ and $r_k$,
such that the Stieltjes transform of any solution $\mu$ to
the indeterminate moment problem is given by
$$
\int_\R \frac{d\mu(x)}{x-z} = - \frac{A(z)\phi(z)-C(z)}
{B(z)\phi(z)-D(z)}
$$
where $\phi$ is any function that is analytic in the
upper half plane having non-negative imaginary part, or
$\phi=\infty$.

Now the solutions to the moment problem corresponding
to $\phi(z) = t\in \R\cup\{\infty\}$ are precisely the
measures that can be obtained from the spectral measure
of a self-adjoint extension of the Jacobi operator $J$.
In particular, this implies that these orthogonality
measures are discrete with support at the zeroes of an entire
function $B(z)t-D(z)$, so that there is no accumulation point.
According to a theorem of M.~Riesz (1923) these are precisely
the measures for which the polynomials are dense in the
corresponding $L^2$-space, cf. Theorem \ref{3340}
for the unbounded case. See Akhiezer \cite{Akhi} and
Simon \cite{Simo} for more
information.
\end{sspar}

\newpage
\section{Doubly infinite Jacobi operators}

\subsection{Doubly infinite Jacobi operators}

\begin{sspar}{\label{410}} We consider an operator on
the Hilbert space $\Hi$, see \ref{210}, of the form
$$
L\, e_k = a_k \, e_{k+1} + b_k\, e_k +
a_{k-1}\, e_{k-1}, \qquad
a_k> 0, \ b_k\in\R,
$$
where $\{ e_k\}_{k\in\Z}$ is the standard orthonormal basis
of $\Hi$ as in \ref{210}. If $a_i=0$ for some $i\in\Z$, then $L$
splits as the direct sum of two Jacobi operators as defined
in \ref{3300}. We call $L$ a Jacobi operator on $\Hi$ or a
doubly infinite Jacobi operator.

The domain ${\mathcal D}(L)$
of $L$ is the dense subspace ${\mathcal D}(\Z)$ of finite linear
combinations of the basis elements $e_k$.  This makes $L$ a
densely defined symmetric operator.
\end{sspar}

\begin{sspar}{\label{420}}
We extend the action of $L$ to an arbitrary vector
$v=\sum_{k=-\infty}^\infty v_ke_k\in\Hi$  by
$$
L^\ast \, v = \sum_{k=-\infty}^\infty (a_k \, v_{k+1}
+ b_k\, v_k + a_{k-1}\, v_{k-1})\, e_k,
$$
which is not an element of $\Hi$ in general.
Define
$$
{\mathcal D}^\ast = \{ v\in\Hi\mid L^\ast v\in \Hi\}.
$$

\begin{lem} $(L^\ast, {\mathcal D}^\ast)$ is the adjoint
of $(L,{\mathcal D}(\Z))$.
\end{lem}

The proof of this Lemma is the same as the proof of
Proposition \ref{3410}.
\end{sspar}

\begin{sspar}{\label{425}}
In particular, $L^\ast$ commutes with complex conjugation,
so its deficiency indices are equal. The solution
space of $L^\ast v=z\, v$ is two-dimensional, since $v$
is completely determined by any initial data $(v_{n-1},v_n)$
for any fixed $n\in\Z$. So the deficiency indices are
equal to $(i,i)$ with $i\in \{ 0,1,2\}$.
{}From \S 2.3 we conclude that $L$ has self-adjoint extensions.
\end{sspar}

\subsection{Relation to Jacobi operators}

\begin{sspar}{\label{430}}
To the operator $L$ we associate two Jacobi operators
$J^+$ and $J^-$ acting on $\Hp$ with orthonormal basis
denoted by $\{f_k\}_{k\in\Zp}$ in order to avoid
confusion. Define
\begin{equation*}
\begin{split}
J^+\, f_k &= \begin{cases}
a_k \, f_{k+1} + b_k\, f_k + a_{k-1}\, f_{k-1},&
\text{for $k\geq 1$,} \\
a_0 \, f_1 + b_0\, f_0, & \text{for $k=0$,}\end{cases} \\
J^-\, f_k &= \begin{cases}
a_{-k-2}\, f_{k+1} + b_{-k-1}\, f_k + a_{-k-1}\, f_{k-1},&
\text{for $k\geq 1$,} \\
a_{-2} \, f_1 + b_{-1}\, f_0, & \text{for $k=0$,}\end{cases}
\end{split}
\end{equation*}
and extend by linearity to ${\mathcal D}(\Zp)$, the space
of finite linear combinations of the basis vectors
$\{f_k\}_{k=0}^\infty$ of $\Hp$. Then $J^{\pm}$ are densely
defined symmetric operators with deficiency indices $(0,0)$
or $(1,1)$ corresponding to whether the associated Hamburger
moment problems is determinate or indeterminate, see \S 3.
The following theorem, due to Masson and Repka \cite{MassR},
relates the deficiency indices of $L$ and $J^\pm$.
\end{sspar}

\begin{thm}\begin{tsspar}\label{435} {\rm (Masson and Repka)}
The deficiency indices of
$L$ are obtained by summing the deficiency indices of $J^+$
and the deficiency indices of $J^-$. \end{tsspar}
\end{thm}

\begin{proof} Let $P_k(z)$, $Q_k(z)$ be generated by the recursion
\begin{equation}
z\, p_k(z) = a_k \, p_{k+1}(z) + b_k\, p_k (z)+
a_{k-1}\, p_{k-1}(z),
\label{eq436}
\end{equation}
subject to the initial conditions
$P_0(z)=1$, $P_{-1}(z)=0$ and $Q_0(z)=0$, $Q_{-1}=1$.
Then, see \ref{425}, any solution of $L^\ast v =z\, v$
can be written as $v_k = v_0 P_k(z) +v_{-1} Q_k(z)$.
Note that $\{ p^+_k\}_{k=0}^\infty$, $p^+_k(z)=P_k(z)$,
are the orthonormal
polynomials corresponding to the Jacobi operator $J^+$ and
that $\{ p^-_k\}_{k=0}^\infty$, $p^-_k(z)=Q_{-k-1}(z)$,
are the orthonormal polynomials corresponding to $J^-$.

Introduce the spaces
\begin{equation}
\begin{split}
S^-_z &= \{ \{f_k\}_{k=-\infty}^\infty \mid
L^\ast f = z\, f\text{\ and\ } \sum_{k=-\infty}^{-1}
|f_k|^2<\infty \}, \\
S^+_z &= \{ \{f_k\}_{k=-\infty}^\infty \mid
L^\ast f = z\, f\text{\ and\ } \sum_{k=0}^\infty
|f_k|^2<\infty \}.
\end{split}
\label{eq409}
\end{equation}
{}From \ref{425} we find $\dim S^\pm_z \leq 2$. The
deficiency space for $L$ is precisely $S^+_z\cap S^-_z$.

{}From the results in \S 3, in particular
Corollary \ref{3100}, Proposition \ref{3350},
Corollary \ref{3410}  and \ref{3430}, we have for
$z\in\C\backslash\R$ that
$\dim S^\pm_z =1$ if and only if $J^\pm$ has
deficiency indices $(0,0)$ and  $\dim S^\pm_z =2$
if and only if $J^\pm$ has deficiency indices $(1,1)$.

Let us now consider the four possible cases.
Case 1: $J^+$ and $J^-$ have deficiency indices
$(1,1)$. Then $\dim S^\pm_z = 2$ and by \ref{425} we
get $S^+_z=S^-_z$, so that $L$ has deficiency indices
$(2,2)$. Case 2: one of $J^\pm$ is essentially self-adjoint.
We can assume that $J^+$ has deficiency indices $(0,0)$ and
that $J^-$ has deficiency indices $(1,1)$. Then
$\dim S^+_z=1$, $\dim S^-_z=2$ and by \ref{245} $S^-_z$ coincides
with the solution space of $L^\ast v=z\, v$, i.e. every
solution of $L^\ast v =z\, v$ is square summable
at $-\infty$.
So (\ref{eq409})
shows that $S^+_z\subset S^-_z$ and hence
$\dim S^+_z\cap S^-_z =1$. So $L$ has deficiency indices
$(1,1)$ in this case.
Case 3: $J^+$ and $J^-$ have deficiency indices
$(0,0)$. Then $\dim S^\pm_z =1$ and we have to show that
$S^+_z\cap S^-_z=\{ 0\}$. Let
$v^\pm=\{v^\pm_k\}_{k=-\infty}^\infty$ span the space $S^\pm_z$,
and we have to show that $v^+$ is not a multiple of $v^-$.
Let $v^\pm_k = C_\pm^P P_k(z) + C_\pm^Q Q_k(z)$.
We calculate
$C^P_\pm$ and $C^Q_\pm$ explicitly in terms of the Stieltjes
transform $w^\pm(z)=\int_\R (x-z)^{-1}\, d\mu^\pm(x)$ of the
spectral measures $\mu^\pm$ for $J^\pm$. Since $J^+$ and $J^-$ are
both essentially self-adjoint, the spaces $S^\pm_z$ are described
in Proposition \ref{3350}. It follows from the recursion
(\ref{eq436}) that the associated polynomials $r^\pm_k(z)$
for $J^\pm$, see \ref{3100}, satisfy
$r^+_k(z)=-a_1^{-1}Q_k(z)$, $k\geq 0$, and
$r^-_k(z) =-a_{-1}^{-1} P_{-k-1}(z)$, $k\geq 0$. By
Proposition \ref{3350}
$v^+_k=f^+_k(z)=w^+(z)P_k(z)-a_1^{-1}Q_k(z)$, $k\geq 0$,
and $v^-_k=f^-_{-k-1}(z)=w^-(z)p^-_{-k-1}(z)+r^-_{-k-1}(z)=
w^-(z)Q_k(z)-a_{-1}^{-1}P_k(z)$. So
$C^P_+=w^+(z)$, $C^Q_+=\frac{-1}{a_1}$,
$C^P_-=\frac{-1}{a_{-1}}$ and $C^Q_-=w^-(z)$,
and consequently
\begin{equation*}
\frac{C_+^P}{C_+^Q} = -a_1\, w^+(z)
={\mathcal O}(\frac{1}{z}), \quad |z|\to\infty, \qquad
\frac{C_-^P}{C_-^Q} = \frac{-1}{a_{-1}w^-(z)}
={\mathcal O}(z), \quad |z|\to\infty.
\end{equation*}
In case $v^+$ is a non-zero multiple of $v^-$ the quotients
 have to be equal.
Since they are also independent of $z$ we find the
required contradiction.
This implies $S^+_z\cap S^-_z=\{ 0\}$ or $L$ has
deficiency indices $(0,0)$.
\end{proof}

\begin{sspar}{\label{440}} As in \ref{3420} we can describe the
sesquilinear form $B$ as introduced in \ref{2120} in terms
of the Wronskian $[u,v]_k = a_k(u_{k+1}v_k-u_kv_{k+1})$. Note
that, as in \S 3, the Wronskian $[u,v]=[u,v]_k$ is independent
of $k$ for $L^\ast u=z\, u$, $L^\ast v=z\, v$, and then
$[u,v]\not=0$ if and only if $u$ and $v$ are linearly
independent solutions. Now, as in \ref{3420},
\begin{equation*}
\begin{split}
&\sum_{k=M}^N (L^\ast u)_k \bar v_k - u_k \overline{(L^\ast v)_k}
\\ &=
\sum_{k=M}^N (a_ku_{k+1}+b_ku_k+a_{k-1}u_{k-1})\bar v_k
-u_k(a_k\bar v_{k+1}+b_k\bar v_k+a_{k-1}\bar v_{k-1}) \\ &=
\sum_{k=M}^N [u,\bar v]_k - [u,\bar v]_{k-1} =
[u,\bar v]_N-[u,\bar v]_{M-1},
\end{split}
\end{equation*}
so that
$$
B(u,v) = \lim_{N\to\infty} [u,\bar v]_N
-\lim_{M\to-\infty}[u,\bar v]_M, \qquad u,v \in {\mathcal D}^\ast.
$$
In particular, if $J^-$ is essentially self-adjoint, we have
$\lim_{M\to-\infty}[u,\bar v]_M=0$, so that in this case the
sesquilinear form $B$ is as in \S 3. In case $J^+$ has
deficiency indices $(1,1)$ we see as in Lemma \ref{3420}
that the same formula for the domains of self-adjoint extensions
of $L$ are valid.

\begin{lem} Let $J^-$, respectively $J^+$, have deficiency
indices $(0,0)$, respectively $(1,1)$, so that $L$ has
deficiency indices $(1,1)$. Then the self-adjoint extensions
of $L$ are given by $(L^\ast, {\mathcal D}_\te)$, $\te\in[0,2\pi)$,
with
$$
{\mathcal D}_\te = \{ v\in {\mathcal D}^\ast \mid
\lim_{N\to\infty} [v, e^{i\te}\Phi_i + e^{-i\te} \Phi_{-i}]_N
=0\}
$$
where $L^\ast \Phi_{\pm i}=\pm i \Phi_{\pm i}$,
$\overline{(\Phi_i)_k} =
(\Phi_{-i})_k$ and
$\| \Phi_{\pm i}\|=1$.
\end{lem}

The proof of Lemma \ref{440} mimics the proof of
Lemma \ref{3420}.
\end{sspar}

\subsection{The Green kernel}

\begin{sspar}{\label{450}} From on we
assume that $J^-$ has deficiency indices $(0,0)$, so that
$J^-$ is essentially self-adjoint and by
Theorem \ref{435} the deficiency indices of $L$ are
$(0,0)$ or $(1,1)$. (The reason for the restriction to this
case is that in case $L$ has deficiency indices $(2,2)$
the restriction of the domain of a self-adjoint extension
of $L$ to the Jacobi operator $J^\pm$ does not in general
correspond to a self-adjoint extension of $J^\pm$, cf.
\cite[Thm.~XII.4.31]{DunfS}.)
Let $z\in\C\backslash\R$ and
choose $\Phi_z\in S^-_z$, so that $\Phi_z$ is determined
up to a constant. We assume $\overline{(\Phi_z)_k} =
(\Phi_{\bar z})_k$, cf. Lemma \ref{440}.
\end{sspar}

\begin{sspar}{\label{460}} Let $\phi_z\in S^+_z$, such that
$\overline{(\phi_z)_k}=(\phi_{\bar z})_k$. We now show that
we may assume
\begin{enumerate}
\item $[\phi_z,\Phi_z]\not= 0$,
\item $\tilde\phi_z$, defined by $(\tilde\phi_z)_k=0$ for
$k<0$ and $(\tilde\phi_z)_k=(\phi_z)_k$ for $k\geq 0$, is
contained in the domain of a self-adjoint extension of $L$.
\end{enumerate}

\vskip0.3truecm
First observe $L^\ast \tilde\phi_z =z\tilde\phi_z
+a_{-1}((\phi_z)_0e_{-1}-(\phi_z)_{-1}e_0)$, so that
$\tilde\phi_z\in{\mathcal D}^\ast$. In case $L$ is essentially
self-adjoint, $(L^\ast, {\mathcal D}^\ast)$ is the unique
self-adjoint extension and (2) is valid. In this case (1)
follows from case 3 in the proof of Theorem \ref{435}.

In case $L$ has deficiency indices $(1,1)$, we have
$S^-_z\subset S^+_z$ and (2) implies (1). Indeed, if (2) holds
and (1) not, then $\Phi_z\in \Hi$ and $\Phi_z=C\phi_z$ is in the
domain of a self-adjoint extension.
Since $L^\ast \Phi_z=z\, \Phi_z$
this shows that the self-adjoint extension would have a non-real
eigenvalue; a contradiction. To show that we can assume
(2) we observe that
$$
\lim_{N\to\infty}
[\tilde\phi_z, e^{i\te}\Phi_i + e^{-i\te} \Phi_{-i}]_N
= \lim_{N\to\infty}
[\phi_z, e^{i\te}\Phi_i + e^{-i\te} \Phi_{-i}]_N,
$$
which exists since $\tilde\phi_z, \Phi_{\pm i}\in {\mathcal D}^\ast$.
If the limit is non-zero, say $K$, we use that
\begin{equation*}
\begin{split}
A(z,\te) &=  \lim_{N\to\infty}
[\Phi_z, e^{i\te}\Phi_i + e^{-i\te} \Phi_{-i}]_N \\ &=
\langle L^\ast \Phi_z, e^{i\te}\Phi_i + e^{-i\te}
\Phi_{-i}\rangle -
\langle \Phi_z, L^\ast( e^{i\te}\Phi_i + e^{-i\te}
\Phi_{-i})\rangle \\
&= e^{-i\te}(z+i)\langle \Phi_z,\Phi_i\rangle +
e^{i\te}(z-i)\langle \Phi_z,\Phi_{-i}\rangle \not= 0.
\end{split}
\end{equation*}
Otherwise, as before, $\Phi_z$ would be in the domain of a
self-adjoint extension of $L$ by Lemma \ref{440}. So that
replacing $\phi_z$ by $\phi_z- \frac{K}{A(z,\te)}\Phi_z$
gives the desired result, since $S^-_z\subset S^+_z$.
\end{sspar}

\begin{sspar}{\label{470}} Let $(L^\prime, {\mathcal D}^\prime)$ be
a self-adjoint extension of $L$, assuming, as before, that
$J^-$ has deficiency indices $(0,0)$.
Let $\phi_z\in S^+_z$, $\Phi_z\in S^-_z$ as in \ref{460}.
We define the Green kernel for $z\in\C\backslash\R$ by
$$
G_{k,l}(z) = \frac{1}{[\phi_z,\Phi_z]}\begin{cases}
(\Phi_z)_k\, (\phi_z)_l, & k\leq l, \\
(\Phi_z)_l\, (\phi_z)_k, & k>l.
\end{cases}
$$
So $\{ G_{k,l}(z)\}_{k=-\infty}^\infty,
\{ G_{k,l}(z)\}_{l=-\infty}^\infty \in \Hi$ and
$\Hi\ni v\mapsto G(z)v$ given by
$$
(G(z)v)_k = \sum_{l=-\infty}^\infty v_l G(z)_{k,l} =
\langle v, \overline{G_{k,\cdot}(z)}\rangle
$$
is well-defined. For $v\in{\mathcal D}(\Z)$ we have
$G(z)v\in {\mathcal D}^\prime$.

\begin{prop} The resolvent of $(L^\prime, {\mathcal
D}^\prime)$ is given by $(L^\prime-z)^{-1}=G(z)$ for
$z\in\C\backslash\R$.
\end{prop}

\begin{proof} Note $\C\backslash\R\subset \rho(L^\prime)$,
because $L^\prime$ is self-adjoint. Hence $(L^\prime-z)^{-1}$
is a bounded operator mapping $\Hi$ onto ${\mathcal D}^\prime$.
So $v\mapsto ((L^\prime-z)^{-1}v)_k$ is a continuous map, hence,
by the Riesz representation theorem,
$((L^\prime-z)^{-1}v)_k = \langle v,
\overline{G_{k,\cdot}(z)}\rangle$ for some $G_{k,\cdot}(z)\in\Hi$.
So it suffices to check $(L^\prime -z )G(z)v=v$ for $v$ in the
dense subspace ${\mathcal D}(\Z)$.
As in the proof of Proposition \ref{3370} we have
\begin{equation*}
\begin{split}
[\phi_z,\Phi_z]&\bigl( (L^\prime -z) G(z)v\bigr)_k =
\sum_{l=-\infty}^{k-1} v_l \bigl( a_k
(\phi_z)_{k+1}+(b_k-z)(\phi_z)_k +a_{k-1}(\phi_z)_{k-1}\bigr)
(\Phi_z)_l \\
& + \sum_{l=k+1}^\infty v_l \bigl( a_k
(\Phi_z)_{k+1}+(b_k-z)(\Phi_z)_k +a_{k-1}(\Phi_z)_{k-1}\bigr)
(\phi_z)_l  \\
& + v_k\bigl( a_k(\Phi_z)_k(\phi_z)_{k+1} +
(b_k-z)(\Phi_z)_k(\phi_z)_k +
a_{k-1}(\Phi_z)_{k-1}(\phi_z)_k\bigr) \\
= & v_k a_k\bigl((\Phi_z)_k(\phi_z)_{k+1} -(\Phi_z)_{k+1}
(\phi_z)_k\bigr) = v_k [\phi_z,\Phi_z]
\end{split}
\end{equation*}
and canceling the Wronskian gives the result.
\end{proof}
\end{sspar}

\begin{sspar}{\label{472}} As in \ref{3380} we can calculate
\begin{equation}
\langle G(z)u, v\rangle = \sum_{k,l=-\infty}^\infty
G_{k,l}(z)u_l\bar v_k =
\frac{1}{[\phi_z,\Phi_z]} \sum_{k\leq l}
(\Phi_z)_k(\phi_z)_l\bigl( u_l\bar v_k+u_k\bar v_l\bigr) (1-\hf
\de_{k,l}),
\label{eq475}
\end{equation}
but in general we cannot pin down the terms of (\ref{eq475})
that will contribute to the spectral measure of
$L^\prime$ in the Stieltjes-Perron inversion formula
of \ref{280}, \ref{2130}. We now consider some examples.
\end{sspar}

\subsection{Example: the Meixner functions}

This subsection is based on Masson and Repka \cite{MassR},
see also Masson \cite{Mass}. We extend to the results of
\cite{MassR} by calculating explicitly the spectral
measure of the Jacobi operator.

\begin{sspar}{\label{480}}
In this example we take for the coefficients of $L$ in \ref{410}
the following;
$$
a_k =a_k(a,\la,\ep) =\sqrt{(\la +k+\ep+1)(k+\ep-\la)},
\qquad b_k = b_k(a,\la,\ep) = 2a(k+\ep).
$$
Without loss of generality we may assume that $0\leq \ep<1$
and $a>0$ by changing to the orthonormal basis $f_k=(-1)^ke_k$
and the operator $-L$. We will do not so in order
to retain symmetry properties.
The conditions $a_k>0$ and $b_k\in\R$ are met if
we require
$a\in\R$, $\ep\in\R$ for $b_k\in\R$ and
either $\la=-\hf+ib$, $b\geq 0$ or $\ep\in[0,\frac{1}{2})$
and $-\hf\leq \la <-\ep$ or $\ep\in(\frac{1}{2},1)$ and
$\la\in(-\frac{1}{2},\ep-1)$ for $a_k>0$. 
Now $L$ is essentially self-adjoint by Theorem
\ref{435}, since $J^\pm$ are
essentially self-adjoint by (i) of Proposition \ref{3450}.
\end{sspar}

\begin{sspar}{\label{490}} In case we choose $\la$, $\ep$ in such a
way that $a_i=0$ for some $i\in\Z$, the corresponding Jacobi
operator on the $\C$-span of $e_k$, $k> i$, can be identified
with the three-term recurrence relation for the
Meixner-Pollaczek, Meixner or Laguerre polynomials depending
on the size of $a$. In case $a_i\not=0$, we can still
consider the corresponding Jacobi operator on $\Hp$ by putting
$a_{-1}=0$. Then the corresponding orthogonal polynomials
are the associated Meixner-Pollaczek, Meixner or Laguerre
polynomials.
\end{sspar}

\begin{sspar}{\label{495}} First note that if we find a
solution $u_k(z)=u_k(z;a,\la,\ep)$ to
\begin{equation}
z\, u_k(z) = a_k(a,\la,\ep)\, u_{k+1}(z) +b_k(a,\la,\ep)\, u_k(z)
+ a_{k-1}(a,\la,\ep)\, u_{k-1}(z)
\label{eq496}
\end{equation}
then $v_k(z) = (-1)^k u_k(-z)$ satisfies
$$
z\, v_k(z) = a_k(-a,\la,\ep)\, v_{k+1}(z) +b_k(-a,\la,\ep)\, v_k(z)
+ a_{k-1}(-a,\la,\ep)\, v_{k-1}(z)
$$
and $w_k(z)= u_{-k}(z;-a,\la,-\ep)$ also satisfies
(\ref{eq496}). Indeed,
introducing a new orthonormal basis $f_k=e_{-k}$ of
$\Hi$ we see that $L$ is given by
$$
L\, f_k=  a_{-k-1}\, f_{k+1} + b_{-k}\, f_k + a_{-k}\, f_{k-1}
$$
and
$$
a_{-k-1}(\la,\ep,a) = a_k(\la,-\ep,-a), \qquad
b_{-k}(\la,\ep,a) = b_k(\la,-\ep,-a).
$$
\end{sspar}

\begin{sspar}{\label{4100}} In order to find explicit solutions to
$L^\ast v =z\, v$ we need the hypergeometric function.
Recall
\begin{equation*}
\begin{split}
{}_2F_1\left( {{a,b}\atop{c}};x\right) &=
\sum_{k=0}^\infty \frac{(a)_k\,(b)_k}{(c)_k\, k!} x^k, \\
(a)_k =a(a+1)\ldots (a+k-1) &= \frac{\Ga(a+k)}{\Ga(a)},
\end{split}
\end{equation*}
where the series is absolutely convergent for $|x|<1$. The
hypergeometric function has a unique analytic continuation
to $\C\backslash[1,\infty)$. It can also be extended non-uniquely
to the cut $[1,\infty)$. See e.g. \cite[Ch.~II]{HTFeen} for all
the necessary material on hypergeometric functions.

\begin{lem} Let $a^2>1$. The functions
\begin{equation*}
\begin{split}
u^{\pm}_k(z) = &
(\pm 2)^{-k} \bigl(\sqrt{a^2-1}\bigr)^{-k}
\frac{\sqrt{\Ga(k+\la+\ep+1)\, \Ga(k+\ep-\la)}}
{\Ga(k+1+\ep\pm z/2\sqrt{a^2-1})} \\ &\times
\, {}_2F_1\left(
{{k+\ep+1+\la,k+\ep-\la}\atop{k+\ep+1\pm z/2\sqrt{a^2-1}}};
\hf\pm\frac{a}{2\sqrt{a^2-1}}\right)
\end{split}
\end{equation*}
and
\begin{equation*}
\begin{split}
v^{\pm}_k(z) = &
(\pm 2)^{k} \bigl(\sqrt{a^2-1}\bigr)^k
\frac{\sqrt{\Ga(-k+\la-\ep+1)\, \Ga(-k-\ep-\la)}}
{\Ga(-k+1-\ep\pm z/2\sqrt{a^2-1})} \\ &\times
\, {}_2F_1\left(
{{-k-\ep+1+\la,-k-\ep-\la}\atop{-k-\ep+1\pm z/2\sqrt{a^2-1}}};
\hf\mp\frac{a}{2\sqrt{a^2-1}}\right)
\end{split}
\end{equation*}
satisfy the recursion {\rm (\ref{eq496})}.
\end{lem}

\begin{proof} The hypergeometric function is a solution to the
hypergeometric differential equation;
$$
x(1-x)y^{\prime\prime} + (c-(a+b+1)x)y^\prime -aby =0.
$$
Now $f(x;a,b,c) =\frac{\Ga(a)\Ga(b)}{\Ga(c)} {}_2F_1(a,b;c;x)$
satisfies $f^\prime(x;a,b,c)=f(x;a+1,b+1,c+1)$. Hence,
\begin{equation*}
\begin{split}
&r_k = f(x;k+1+\ep+\la,k+\ep-\la,k+1+\ep-y), \\
&x(1-x)r_{k+1} + (k+\ep-y-2(k+\ep)x)r_k -
(k+\ep+\la)(k+\ep-\la-1)r_{k-1}=0.
\end{split}
\end{equation*}
Replace $x=\hf - \frac{a}{2\sqrt{a^2-1}}$ and $y=z/2\sqrt{a^2-1}$.
Let
$$
r_k=(-2)^k \bigl(a^2-1)^{\hf k}
\sqrt{\Ga(k+\ep+1+\la)\Ga(k+\ep-\la)}p_k
$$
then $p_k$ satisfies (\ref{eq496}).

This proves the lemma for $u^-_k(z)$. The case
$u^+_k(z)$ follows
by replacing $a$ by $-a$, and applying \ref{495}.
Replacing $k$ by $-k$, $a$ by $-a$ and $\ep$ by $-\ep$
and using \ref{495} gives the other sets of
solutions to (\ref{eq496}).
\end{proof}

{}From now on we assume that
$a^2>1$ in order not to complicate matters. The case $a^2=1$
corresponds to the Laguerre functions, and the
case $0\leq a^2<1$ corresponds to the Meixner-Pollaczek
functions, see \cite{MassR} and \ref{4160}.
\end{sspar}

\begin{sspar}{\label{4110}} In order to find the asymptotic
behaviour of the solutions in Lemma \ref{4100} we first
use
$$
{}_2F_1\left( {{a,b}\atop{c}};x\right) = (1-x)^{c-a-b}
{}_2F_1\left( {{c-a,c-b}\atop{c}};x\right),
$$
so that
\begin{equation*}
\begin{split}
 {}_2F_1&\left(
{{k+\ep+1+\la,k+\ep-\la}\atop{k+\ep+1\pm z/2\sqrt{a^2-1}}};
\hf\pm\frac{a}{2\sqrt{a^2-1}}\right)  \\  &=
(\hf\mp \frac{a}{2\sqrt{a^2-1}})^{-k-\ep\pm z/2\sqrt{a^2-1}}
\bigl( 1+{\mathcal
O}(\frac{1}{k})\bigr), \qquad k\to\infty.
\end{split}
\end{equation*}
Use
$$
\frac{\Ga(z+a)}{\Ga(z+b)} = z^{a-b}\bigl(1+{\mathcal
O}(\frac{1}{z})\bigr), \qquad z\to\infty,\ |\arg(z)|<\pi,
$$
to find
$$
\frac{\sqrt{\Ga(k+\la+\ep+1)\, \Ga(k+\ep-\la)}}{|\Ga(k+1+\ep+y)|} =
k^{-\Re y-\hf}\bigl( 1+ {\mathcal O}(\frac{1}{k})\Bigr),
\qquad k\to\infty,
$$
so that we find the asymptotic behaviour
$$
|u^\pm_k(z)| =
(\hf\mp \frac{a}{2\sqrt{a^2-1}})^{-\ep\pm \Re z/2\sqrt{a^2-1}}
k^{-\hf \mp \Re z/2\sqrt{a^2-1}}
\, |-a\pm \sqrt{a^2-1}|^{-k}
\bigl( 1+{\mathcal
O}(\frac{1}{k})\bigr),
$$
as $k\to\infty$.
This implies that $u^+(z)\in S^+_z$ if $a<-1$ and
$u^-(z)\in S^+_z$ if $a>1$.

Similarly we obtain the asymptotic behaviour
$$
|v^\pm_{-k}(z)| =
(\hf\pm \frac{a}{2\sqrt{a^2-1}})^{\ep\pm \Re z/2\sqrt{a^2-1}}
k^{-\hf \mp \Re z/2\sqrt{a^2-1}}
(a\pm \sqrt{a^2-1})^{-k}
\bigl( 1+{\mathcal
O}(\frac{1}{k})\bigr),
$$
as $k\to\infty$,
so that $v^-(z)\in S^-_z$ if $a<-1$ and $v^+(z)\in S^-_z$
if $a>1$.

Note that we have completely determined $S^\pm_z$,
hence $\phi_z$ and $\Phi_z$, since these
spaces are one-dimensional.
\end{sspar}

\begin{sspar}{\label{4115}}
If we reparametrise the parameter
$a=\hf(s+s^{-1})$, then
$$
 \{ s, s^{-1}\} =
\{ a+\sqrt{a^2-1},a-\sqrt{a^2-1}\}.
$$
Note that $a^2-1=\frac{1}{4}(s-s^{-1})^2$.
In this case we can let $u^\pm(z)$ and $v^\pm(z)$ correspond
to hypergeometric series with only $s^{\pm 1}$-dependence.
The results can be written somewhat nicer after
transforming Lemma \ref{4100} by
$$
{}_2F_1\left( {{a,b}\atop{c}};x\right) = (1-x)^{-a}
{}_2F_1\left( {{a,c-b}\atop{c}};\frac{x}{x-1}\right).
$$
We leave this to the reader.
\end{sspar}

\begin{sspar}{\label{4120}} Note that $v^\pm(z)$ are linearly
independent solutions to (\ref{eq496}), because they display
different asymptotic behaviour. Since the solution space
to (\ref{eq496}) is two-dimensional, we see that there
exist constants such that
$$
u^+_k(z) = A^+(z) v^+_k(z) + B^+(z) v^-_k(z), \qquad
u^-_k(z) = A^-(z) v^+_k(z) + B^-(z) v^-_k(z).
$$
Or stated differently, relations between hypergeometric series
of argument $x$ and $1-x$. Relations connecting
hypergeometric series are very classical.
We use
\begin{equation*}
\begin{split}
&{}_2F_1\left( {{a,b}\atop{a+b+1-c}};1- x\right) =
A  {}_2F_1\left( {{a,b}\atop{c}}; x\right)
+ B  x^{1-c} (1-x)^{c-a-b}
{}_2F_1\left( {{1-a,1-b}\atop{2-c}}; x\right),  \\
&A = \frac{\Ga(a+b+1-c)\Ga(1-c)}{\Ga(a+1-c)\Ga(b+1-c)}, \qquad
B = \frac{\Ga(a+b+1-c)\Ga(c-1)}{\Ga(a)\, \Ga(b)}
\end{split}
\end{equation*}
with $a\mapsto -k-\ep+1+\la$, $b\mapsto -k-\ep-\la$,
$c\mapsto -k-\ep+1+y$, $x\mapsto \hf -\frac{a}{2\sqrt{a^2-1}}$.
Then the first ${}_2F_1$ is as in $v^-_k(z)$, the second one
as in $v^+_k(z)$, and the third one as in $u^-_k(z)$.
This gives, using $y=z/2\sqrt{a^2-1}$,
\begin{equation*}
\begin{split}
B^-(z) =& (-2)^{-k} \bigl(\sqrt{a^2-1}\bigr)^{-k}
\frac{\Ga(-k+1-\ep-y)}{\sqrt{\Ga(-k+\la-\ep+1)\Ga(-k-\ep-\la)}}
\\ &\hskip-42pt
\rightline{\text{from factor of $v_k^-(z)$}} \\
\times&
\frac{\Ga(-k-\ep+1+\la)\Ga(-k-\ep-\la)}
{\Ga(1-k-\ep-y)\Ga(-k-\ep+y)}
\bigl(\hf-\frac{a}{2\sqrt{a^2-1}}\bigr)^{y-k-\ep}
\bigl(\hf+\frac{a}{2\sqrt{a^2-1}}\bigr)^{-y-k-\ep}
\\ &\hskip-42pt \rightline{\text{from $B$}} \\
\times&(-2)^{-k} \bigl(\sqrt{a^2-1}\bigr)^{-k}
\frac{\sqrt{\Ga(k+\la+\ep+1)\Ga(k+\ep-\la)}}{\Ga(k+1+\ep-y)}
\\ & \hskip-42pt\rightline{\text{from factor of $u_k^-(z)$.}}
\end{split}
\end{equation*}
Now $B^-(z)$ has to be independent of $k$. After canceling
common factors we use
$\Ga(z)\Ga(1-z)=\pi/\sin \pi z$, $z\not\in\Z$, to obtain
$$
B^-(z) =
\Bigl( \frac{\sqrt{a^2-1}-a}{\sqrt{a^2-1}+a}\Bigr)^y
(4(1-a^2))^\ep
\frac{\sin (\pi(y-\ep))}{\sqrt{\sin(\pi(\ep-\la))\,
\sin(\pi(-\ep-\la))}}.
$$
Similarly, also using $\Ga(z+1)=z\Ga(z)$, we find
$$
A^-(z) = \frac{\pi}{\Ga(1+\la-y)\Ga(-\la-y)}
 \Bigl( \frac{\sqrt{a^2-1}-a}{\sqrt{a^2-1}+a}\Bigr)^y
\frac{(4(1-a^2))^\ep}{\sqrt{\sin(\pi(\ep-\la))\,
\sin(\pi(-\ep-\la))}}.
$$
\end{sspar}

\begin{sspar}{\label{4130}} Next we calculate the Wronskians
$[v^-(z),v^+(z)]$ and $[u^-(z),v^+(z)]$.
To calculate the first Wronskian
observe that we can take the limit $k\to\infty$
in
$$
[v^-(z),v^+(z)] = a_{-k}\bigl( v^-_{-k+1}(z)v^+_{-k}(z) -
v^-_{-k}(z)v^+_{-k+1}(z)\bigr)
$$
using the asymptotic behaviour of \ref{4110}. This gives,
using $y=z/2\sqrt{a^2-1}$,
\begin{equation*}
\begin{split}
[v^-(z),v^+(z)] = & a_{-k}
(\hf- \frac{a}{2\sqrt{a^2-1}})^{\ep-y}
(\hf+ \frac{a}{2\sqrt{a^2-1}})^{\ep+y}
\\ &\times
\Bigl( (k-1)^{y-\hf} (a-\sqrt{a^2-1})^{1-k} k^{-\hf-y}
(a+\sqrt{a^2-1})^{-k} \\ &\qquad -
k^{y-\hf} (a-\sqrt{a^2-1})^{-k} (k-1)^{-\hf-y}
(a+\sqrt{a^2-1})^{1-k}\Bigr)
\end{split}
\end{equation*}
and taking out common factors and using that $a_{-k}=k(1+{\mathcal
O}(\frac{1}{k}))$ gives
$$
[v^-(z),v^+(z)] = -2\sqrt{a^2-1} (4(1-a^2))^{-\ep}
\Bigl( \frac{\sqrt{a^2-1}+a}{\sqrt{a^2-1}-a}\Bigr)^y.
$$

{}From \ref{4120} it follows that, $y=z/2\sqrt{a^2-1}$,
\begin{equation*}
[u^-(z),v^+(z)]
=  B^-(z) [v^-(z),v^+(z)]
= -2\sqrt{a^2-1}
\frac{\sin (\pi(y-\ep))}{\sqrt{\sin(\pi(\ep-\la))\,
\sin(\pi(-\ep-\la))}}.
\end{equation*}
We see that the Wronskian $[u^-(z),v^+(z)]$, as a function of
$z$, has no poles and it has zeroes at
$z=2(\ep+l)\sqrt{a^2-1}$, $l\in\Z$.
\end{sspar}

\begin{sspar}{\label{4140}} With all these preparations we
can calculate the spectral measure for the doubly
infinite Jacobi operator $L$ with coefficients as in
\ref{480}, where we assume $a>1$. So we can take
$\phi_z=u^-(z)$ and $\Phi_z=v^+(z)$. And we see that
these solutions are analytic in $z$, since
$\Ga(c)^{-1} \, {}_2F_1(a,b;c;z)$ is analytic in $c$.
Hence, the only contribution in the spectral measure,
cf. \ref{472}, comes from the zeroes of the Wronskian,
so that the spectrum is purely discrete.

\begin{thm} The operator
$L$ on $\Hi$ defined by {\rm \ref{410}}, {\rm \ref{480}}
is essentially self-adjoint. For $a>1$ its
self-adjoint extension has completely discrete spectrum
$\{ 2(\ep+l)\sqrt{a^2-1}\}_{l\in\Z}$, and in particular
the set
$$
\bigl\{ u^-(2(\ep+l)\sqrt{a^2-1}) \mid l\in\Z\bigr\}
$$
constitutes a complete orthogonal basis of $\Hi$. Moreover,
\begin{equation*}
\begin{split}
&\|  u^-(2(\ep+l)\sqrt{a^2-1})\|^2 = \\
&\Bigl(
\frac{4(a^2-1)(a-\sqrt{a^2-1})}{a+\sqrt{a^2-1}}\Bigr)^{\ep}
\Bigl( \frac{a-\sqrt{a^2-1}}{a+\sqrt{a^2-1}}\Bigr)^l
\Ga(\ep+l-\la)\Ga(1+\ep+l+\la).
\end{split}
\end{equation*}
\end{thm}

\noindent
{\it Remark} (i) Note that
hypergeometric expression for $u^-$ as
in Lemma \ref{4100} displays Bessel coefficient behaviour,
and we can think of the orthogonality relations as
a natural extension
of the Hansen-Lommel ortho\-go\-na\-lity
re\-la\-tions $\sum_{k=-\infty}^\infty
J_{k+n}(z)J_k(z)=\de_{n,0}$ for the Bessel functions,
see \ref{4150}.
\par\noindent
(ii) Note that the spectrum is independent of $\la$.

\begin{proof} Put $\phi_z=u^-(z)$ and $\Phi_z=v^+(z)$, then
we see that, see (\ref{eq475}),
$$
\langle G(z)u, v\rangle = \sum_{k,l=-\infty}^\infty
G_{k,l}(z)u_l\bar v_k =
\frac{1}{[\phi_z,\Phi_z]} \sum_{k\leq l}
(\Phi_z)_k(\phi_z)_l\bigl( u_l\bar v_k+u_k\bar v_l\bigr) (1-\hf
\de_{k,l}).
$$
Hence, the only contribution in
$$
E_{u,v}\bigl( (a,b)\bigr) = \lim_{\de\downarrow 0}
\lim_{\ep\downarrow 0} \frac{1}{2\pi i}
\int_{a+\de}^{b-\de}
\langle G(x+i\ep)u,v\rangle - \langle G(x-i\ep)u,v\rangle \, dx
$$
can come from the zeroes of the Wronskian. Put
$x_l=2(\ep+l)\sqrt{a^2-1}$, then for $(a,b)$ containing
precisely one of the $x_l$'s we find
$$
E_{u,v}\bigl( (a,b)\bigr) =
E_{u,v}\bigl( \{x_l\} \bigr) = - \frac{1}{2\pi i}\oint_{(x_l)}
\langle G(z)u,v\rangle \, dz.
$$
Note that the minus sign comes from the clockwise
orientation of the rectangle with upper side $x+i\ep$, $x\in[a,b]$,
and lower side $x-i\ep$, $x\in[a,b]$.
This residue can be calculated and we find
\begin{equation*}
\begin{split}
\text{Res}_{z=x_l} \frac{1}{[\phi_z,\Phi_z]} &=
\frac{-\sqrt{\sin(\pi(\ep-\la))\, \sin(\pi(-\ep-\la))}}
{2\sqrt{a^2-1}}
\text{Res}_{z=x_l} \frac{1}{\sin(-\ep\pi+\pi z/2\sqrt{a^2-1})}\\
&= \frac{(-1)^{l+1}}{\pi}
\sqrt{\sin(\pi(\ep-\la))\, \sin(\pi(-\ep-\la))}
\end{split}
\end{equation*}

Since $u_k^-(x_l)=A^-(x_l)v^+_k(x_l)$, because the zeroes
of the Wronskian $[u^-(z),v^+(z)]$ correspond to the
zeroes of $B^-(z)$, see \ref{4130}, we see that
$\phi_z$ and $\Phi_z$ are multiples of each other for $z=x_l$.
So we can symmetrise the sum in $\langle G(z)u,v\rangle$
again, and we find
$$
E_{u,v}\bigl( \{x_l\} \bigr) =
\frac{(-1)^l}{\pi A^-(x_l)}
\sqrt{\sin(\pi(\ep-\la))\, \sin(\pi(-\ep-\la))}
\langle u, \phi_{x_l}\rangle \langle v,\phi_{x_l}\rangle .
$$
Using the explicit expression for $A^-(x_l)$ of \ref{4120}
and the reflection identity for the $\Ga$-function we arrive at
$$
E_{u,v}\bigl( \{x_l\} \bigr) =
\Bigl( \frac{4(a^2-1)(a-\sqrt{a^2-1})}{a+\sqrt{a^2-1}}\Bigr)^{-\ep}
\Bigl( \frac{a+\sqrt{a^2-1}}{a-\sqrt{a^2-1}}\Bigr)^l
\frac{\langle u, \phi_{x_l}\rangle \langle v,\phi_{x_l}\rangle }
{\Ga(\ep+l-\la)\Ga(1+\ep+l+\la)}.
$$

Note that in particular, $\phi_{x_l}$ are eigenvectors of
$L$ for the eigenvalue $x_l$, and moreover these eigenspaces are
$1$-dimensional. Hence
$E(\{x_l\})\phi_{x_l}=\phi_{x_l}$.
Recall $E_{u,v}\bigl( \{x_l\} \bigr) =
\langle E\bigl( \{x_l\} \bigr) u,v\rangle$, and take
$u=v=\phi_{x_l}$ to find
$$
\| \phi_{x_l}\|^2 =
\Bigl( \frac{4(a^2-1)(a-\sqrt{a^2-1})}{a+\sqrt{a^2-1}}\Bigr)^{-\ep}
\Bigl( \frac{a+\sqrt{a^2-1}}{a-\sqrt{a^2-1}}\Bigr)^l
\frac{\|\phi_{x_l}\|^4 }
{\Ga(\ep+l-\la)\Ga(1+\ep+l+\la)}
$$
from which the squared norm follows.
Using $E(\{x_l\})E(\{ x_m\}) =\de_{l,m} E(\{x_l\})$
and self-adjointness of $E$ gives
\begin{equation*}
\begin{split}
\langle \phi_{x_l},\phi_{x_m}\rangle &=
\langle E(\{x_l\})\phi_{x_l},E(\{ x_m\})\phi_{x_m}\rangle =
\langle E(\{ x_m\})E(\{x_l\})\phi_{x_l},\phi_{x_m}\rangle \\ &=
\de_{l,m} \langle E(\{x_l\})\phi_{x_l},\phi_{x_m}\rangle =
\de_{l,m}\langle \phi_{x_l},\phi_{x_m}\rangle
\end{split}
\end{equation*}
which proves the orthogonality.
\end{proof}
\end{sspar}

\begin{sspar}{\label{4150}} The orthogonality relations
arising from Theorem \ref{4140} can be worked out and they
give
\begin{equation*}
\begin{split}
\sum_{k=-\infty}^\infty &\frac{\Ga(k+\la+\ep+1)\Ga(k+\ep-\la)}
{(4(a^2-1))^k} \\
&\times \frac{1}{\Ga(k+1-l)} {}_2F_1\left(
{{k+\ep+\la+1,k+\ep-\la}\atop{k+1-l}};
\hf-\frac{a}{2\sqrt{a^2-1}}\right)\\
&\times
\frac{1}{\Ga(k+1-p)} {}_2F_1\left(
{{k+\ep+\la+1,k+\ep-\la}\atop{k+1-p}};
\hf-\frac{a}{2\sqrt{a^2-1}}\right) \\
&= \de_{p,l}
\Bigl( \frac{4(a^2-1)(a-\sqrt{a^2-1})}{a+\sqrt{a^2-1}}\Bigr)^{\ep}
\Bigl( \frac{a-\sqrt{a^2-1}}{a+\sqrt{a^2-1}}\Bigr)^l
\Ga(\ep+l-\la)\Ga(1+\ep+l+\la).
\end{split}
\end{equation*}
Since Theorem \ref{4140} also gives the completeness of this
set of vectors, the dual orthogonality relations also hold,
i.e.
\begin{equation*}
\begin{split}
\sum_{l=-\infty}^\infty&
\Bigl( \frac{a+\sqrt{a^2-1}}{a-\sqrt{a^2-1}}\Bigr)^l
\frac{1}{\Ga(\ep+l-\la)\Ga(1+\ep+l+\la)} \\
&\times \frac{1}{\Ga(k+1-l)} {}_2F_1\left(
{{k+\ep+\la+1,k+\ep-\la}\atop{k+1-l}};
\hf-\frac{a}{2\sqrt{a^2-1}}\right) \\
&\times \frac{1}{\Ga(m+1-l)} {}_2F_1\left(
{{m+\ep+\la+1,m+\ep-\la}\atop{m+1-l}};
\hf-\frac{a}{2\sqrt{a^2-1}}\right) \\
& = \de_{k,m} \frac{(4(a^2-1))^k}{\Ga(k+\la+\ep+1)\Ga(k+\ep-\la)}
\Bigl( \frac{4(a^2-1)(a-\sqrt{a^2-1})}{a+\sqrt{a^2-1}}\Bigr)^{\ep}.
\end{split}
\end{equation*}
\end{sspar}

\vskip.1truecm
\begin{sspar}{\label{4160}} {\it Exercise.} Use the results
of \cite{MassR} in order to calculate explicitly the
spectral measures for the case $|a|\leq 1$. In this case
the spaces $S^\pm_z$ are spanned by different solutions
depending on the sign of $\Im z$.

In case $|a|\leq 1$ the situation changes considerably,
and let us state briefly some of the results needed
for the Meizner-Pollaczek case, i.e. $|a|<1$. The
analogue of Lemma \ref{4100} is proved in entirely the same
way, and we find that
\begin{equation*}
\begin{split}
U^{\pm}_k(z) &= (\pm 1)^k (2i\sin\psi)^{-k}
\frac{\sqrt{\Ga(k+1+\la+\ep)\Ga(k+\ep-\la)}}
{\Ga(k+1+\ep\mp iz)}\\ &\qquad\times
\, {}_2F_1\left( {{k+1+\la+\ep, k+\ep-\la}\atop{k+1+\ep\mp iz}};
\frac{1}{1-e^{\pm 2i\psi}}\right), \\
V^{\pm}_k(z) &= (\pm 1)^k (2i\sin\psi)^k
\frac{\sqrt{\Ga(1-k+\la-\ep)\Ga(-k-\ep-\la)}}
{\Ga(1-k-\ep\mp iz)}\\ &\qquad\times
\, {}_2F_1\left( {{1-k+\la-\ep, -k-\ep-\la}\atop{1-k-\ep\mp iz}};
\frac{1}{1-e^{\mp 2i\psi}}\right),
\end{split}
\end{equation*}
satisfy the recurrence relation
$$
(2\sin\psi)\, z\, u_k(z) = a_k(\cos\psi,\la,\ep)\,
u_{k+1}(z) + b_k(\cos\psi,\la,\ep)\, u_k(z)
+ a_{k-1}(\cos\psi,\la,\ep)\,
u_{k-1}(z)
$$
with the explicit
values for $a_k$ and $b_k$ of \ref{480} with $a=\cos\psi$,
$0<\psi<\pi$. This gives four linearly independent solutions
of the recurrence relation. The connection formulae for
these four solutions are given by
\begin{equation*}
U^+_k(z) = A^+(z)\, V_k^+(z) + B^+(z)\, V_k^-(z), \qquad
U^-_k(z) = A^-(z)\, V_k^+(z) + B^-(z)\, V_k^-(z)
\end{equation*}
and it follows easily from the explicit expressions 
and the assumptions on $\la$ and $\ep$ 
of \ref{480} that
$\overline{U_k^+(z)}=U_k^-(\bar z)$ and
$\overline{V_k^+(z)}=V_k^-(\bar z)$, implying that
$A^-(z)=\overline{B^+(\bar z)}$ and
$B^-(z)=\overline{A^+(\bar z)}$. The same formula
for hypergeometric series as in \ref{4120}
can be used to find the connection coefficients;
\begin{equation*}
\begin{split}
A^+(z) & = (2\sin\psi)^{2\ep}e^{2z(\psi-\frac{\pi}{2})}
\frac{\sqrt{\Ga(-\ep-\la)\Ga(1+\la-\ep)\Ga(1+\la+\ep)\Ga(\ep-\la)}}
{\Ga(iz-\ep)\Ga(1+\ep-iz)}, \\
B^+(z) & = (2\sin\psi)^{2\ep}e^{2z(\psi-\frac{\pi}{2})}
\frac{\sqrt{\Ga(-\ep-\la)\Ga(1+\la-\ep)\Ga(1+\la+\ep)\Ga(\ep-\la)}}
{\Ga(\la+1-iz)\Ga(-\la-iz)}.
\end{split}
\end{equation*}

The asymptotic behaviour can be determined
as in \ref{4110}, and we find that $\phi_z=U^+(z)$ for
$\Im z>0$ and $\phi_z=U^-(z)$ for $\Im z<0$ and
$\Phi_z=V^+(z)$ for $\Im z>0$ and $\Phi_z=V^-(z)$ for
$\Im z<0$. The explicit asymptotic behaviour can be used as
in \ref{4130} to find the Wronskian
$$
[V^-(z),V^+(z)] = -i\, (2\sin\psi)^{1-2\ep}
e^{-2z(\psi-\frac{\pi}{2})},
$$
and from this expression and the connection coefficients
we can calculate all Wronskians needed. In order to find
the spectral measure we have to investigate
the limit $\ep\downarrow 0$ of
$\langle G(x+i\ep)u,v\rangle - \langle G(x-i\ep)u,v\rangle$,
$x\in\R$, and for
this we consider, using the connection formulae, 
\begin{equation*}
\begin{split}
&\qquad\qquad\qquad\qquad
\frac{V^+_k(x)U^+_l(x)}{B^+(x)[V^-(x),V^+(x)]}
- \frac{V^-_k(x)U^-_l(x)}{A^-(x)[V^+(x),V^-(x)]} = \\
& \frac{V^+_k(x)V^-_l(x)+V^-_k(x)V^+_l(x)}{[V^-(x),V^+(x)]}
+ \frac{A^+(x)A^-(x)V^+_l(x)V^+_k(x) +
B^-(x)B^+(x)V^-_k(x)V^-_l(x)}
{A^-(x)B^+(x)[V^-(x),V^+(x)]},
\end{split}
\end{equation*}
which is obviously symmetric in $k$ and $l$. Hence we
can antisymmetrise the sum for the spectral measure
and we find for $u,v\in\Hi$ 
\begin{equation*}
\begin{split}
2\pi i\langle u,v\rangle = \int_\R &
\Bigl( A^-(x)B^+(x)\langle u, V^+(x)\rangle
\langle V^+(x),v\rangle +
A^-(x)B^+(x)\langle u, V^-(x)\rangle
\langle V^-(x),v\rangle\\ &+
A^+(x)A^-(x)\langle u, V^-(x)\rangle
\langle V^+(x),v\rangle +
B^+(x)B^-(x)\langle u, V^+(x)\rangle
\langle V^-(x),v\rangle\Bigr)
\\ &\qquad \times\frac{dx}{A^-(x)B^+(x)[V^-(x),V^+(x)]} \\
= \int_\R & \Bigl( \langle u, U^-(x)\rangle\langle U^-(x), v\rangle
+(A^-(x)B^+(x)-A^+(x)B^-(x))
\langle u, V^-(x)\rangle\langle V^-(x), v\rangle\Bigr)
\\ &\qquad \times\frac{dx}{A^-(x)B^+(x)[V^-(x),V^+(x)]}
\end{split}
\end{equation*}
describing the spectral measure, where we have used the
relations between $U^\pm_k(z)$ and $V^\pm_k(x)$ for the
second equality.  
So we see that the spectrum
of the corresponding operator is $\R$, and by inserting
the values for the Wronskian and the connection coefficients
we see that the spectral
measure is described by the following integral;
\begin{equation*}
\begin{split}
&\langle u,v\rangle = \\ \frac{1}{2\pi}
\int_\R &\frac{|\Ga(\la+1-ix)\Ga(-\la-ix)|^2}
{\Ga(-\ep-\la)\Ga(1+\la-\ep)\Ga(1+\la+\ep)\Ga(\ep-\la)}
(2\sin\psi)^{-1-2\ep}e^{-2x(\psi-\frac{\pi}{2})}
\langle u,U^-(x)\rangle \langle U^-(x),v\rangle \\ & +
\Bigl( 1-\frac{|\Ga(\la+1-ix)\Ga(-\la-ix)|^2}
{|\Ga(ix-\ep)\Ga(1+\ep-ix)|^2}\Bigr)
(2\sin\psi)^{-1+2\ep}e^{2x(\psi-\frac{\pi}{2})}
\langle u,V^-(x)\rangle \langle V^-(x),v\rangle\, dx.
\end{split}
\end{equation*}
We see that the spectral projection is on a two-dimensional
space of generalised eigenvectors. Note that the general
theory ensures the positivity of the measure in case
$u=v$, and we
can also check directly that, under the conditions on
$\la$ and $\ep$ as in \ref{480}, the second term in the
integrand is indeed positive.
\end{sspar}

\subsection{Example: the basic hypergeometric
difference equation}

This example is based on Appendix A in \cite{KoelSsu}, which
was greatly motivated by Kakehi \cite{Kake} and unpublished notes
by Koornwinder. The transform described in this section
has been obtained from its quantum $SU(1,1)$ group theoretic
interpretation, see \cite{Kake}, \cite{KoelSsu} for references.
On a formal level the result can be obtained as a limit case
of the orthogonality of the Askey-Wilson polynomials, see
\cite{KoelSNATO} for a precise formulation.
The limit transition descibed in \ref{4250} is motivated
from the fact that the Jacobi operators in this example and
the previous example play the same role.

\begin{sspar}{\label{4170}} We take the coefficients as
$$
a_k = \hf \sqrt{
(1-\frac{q^{-k}}{r})(1-\frac{cq^{-k}}{d^2r})}, \qquad
b_k =\frac{q^{-k}(c+q)}{2dr},
$$
where we assume $0<q<1$, and $r<0$, $c>0$, $d\in\R$.
This assumption is made in order to get the expression under
the square root sign positive. There are more possible choices
in order to achieve this, see \cite[App.~A]{KoelSsu}.
Note that $a_k$ and $b_k$ are bounded for $k<0$, so that
$J^-$ is self-adjoint. Hence, the deficiency indices of $L$ are
$(0,0)$ or $(1,1)$ by Theorem \ref{435}.
\end{sspar}

\begin{sspar}{\label{4175}} In order to write down solutions of
$Lu=z\, u$ we need the basic hypergeometric series.
Define
$$
(a;q)_k = \prod_{i=0}^{k-1}(1-aq^i), \qquad k\in\Zp\cup\{\infty\},
\qquad
(a_1,\ldots,a_n;q)_k = \prod_{j=1}^n (a_j;q)_k,
$$
and the basic hypergeometric series
$$
{}_2\vp_1\left( {{a,b}\atop{c}};q,x\right) =
\sum_{k=0}^\infty \frac{(a,b;q)_k}{(q,c;q)_k} x^k.
$$
The radius of convergence is $1$, and there exists a unique
analytic continuation to $\C\backslash[1,\infty)$. See
Gasper and Rahman \cite{GaspR} for all the necessary information
on basic hypergeometric series.

\begin{lem} Put
\begin{equation*}
\begin{split}
w_k^2 &= d^{2k} \frac{(cq^{1-k}/d^2r;q)_\infty}
{(q^{1-k}/r;q)_\infty}, \\
f_k(\mu(y)) &=
{}_2\vp_1\left( {{dy,d/y}\atop c};q,rq^k\right),
\qquad c\not\in q^{-\Zp},\quad \mu(y)=\hf(y+y^{-1}), \\
g_k(\mu(y)) &= q^kc^{-k} \, {}_2\vp_1\left(
{{qdy/c, qd/cy}\atop{q^2/c}};q, rq^k\right),
\quad \mu(y)=\hf(y+y^{-1}), \\
F_k(y) &= (dy)^{-k} \, {}_2\vp_1\left({{dy,qdy/c}\atop{qy^2}};q,
\frac{q^{1-k}c}{d^2r}\right), \qquad y^2\not\in q^{-\N},
\end{split}
\end{equation*}
then, with $z=\mu(y)$, we have that
$u_k(z)= w_kf_k(\mu(y))$, $u_k(z)= w_kg_k(\mu(y))$,
$u_k(z)=w_k F_k(y)$ and
$u_k(z)=w_kF_k(y^{-1})$
define solutions to
$$
z\, u_k(z) = a_k \, u_{k+1}(z) + b_k \, u_k(z) +
a_{k-1}\, u_{k-1}(z).
$$
\end{lem}

\begin{proof} Put $u_k(z)=w_kv_k(z)$, then $v_k(z)$ satisfies
$$
2z\, v_k(z) = (d-\frac{cq^{-k}}{dr})\, v_{k+1}(z)
+q^{-k}\frac{c+q}{dr}\, v_k(z)
+ (d^{-1}-\frac{q^{1-k}}{dr})\, v_{k-1}(z)
$$
and this is precisely the second order $q$-difference equation
that has the solutions given, see \cite[exerc.~1.13]{GaspR}.
\end{proof}
\end{sspar}

\begin{sspar}{\label{4180}} The asymptotics of the solutions
of Lemma \ref{4175} can be given as follows.
First observe that $w_{-k} =d^{-k}$ as $k\to\infty$, and
using
$$
 w_k^2=c^k
\frac{(r q^k, d^2r/c, cq/d^2r;q)_\infty}
{(d^2rq^k/c, r, q/ r;q)_\infty} \Rightarrow
w_k={\mathcal O}(c^{\hf k}), \ k\to\infty.
$$
Now $f_k(\mu(y))= {\mathcal O}(1)$ as $k\to\infty$, and
$g_k(\mu(y)) = {\mathcal O}((q/c)^k)$ as $k\to \infty$.
Similarly, $F_{-k}(y)={\mathcal O}((dy)^k)$ as $k\to\infty$.

\begin{prop} The operator $L$ is essentially
self-adjoint for $0<c\leq q^2$, and $L$ has deficiency
indices $(1,1)$ for $q^2<c<1$. Moreover,
for $z\in\C\backslash\R$ the one-dimensional
space $S^-_z$ is spanned by $wF(y)$ with $\mu(y)=z$ and $|y|<1$.
For $0<c\leq q^2$ the one-dimensional
space $S^+_z$ is spanned by $wf(z)$, and for $q^2<c<1$
the two-dimensional
space $S^+_z$ is spanned by $wf(z)$ and $wg(z)$.
\end{prop}

\begin{proof} In \ref{4170} we have already observed that the
deficiency indices of $L$ are $(0,0)$ or $(1,1)$. Now
$2a_k=q^{-k}\sqrt{c/d^2r^2}-\hf(r+d^2r/c)+{\mathcal O}(q^k)$,
$k\to\infty$, shows that
the boundedness condition of Proposition \ref{3450}(ii)
is satisfied if the
coefficient of $q^{-k}$ in $a_k + a_{k-1}\pm b_k$ is non-positive.
Since $c>0$, $dr\in\R$, this is
the case when $(1+q)\sqrt{c}\leq c+q$.
For $0<c\leq q^2$ the inequality holds, so that by
Proposition \ref{3450}(ii) also $J^+$, and hence
$L$ by Theorem \ref{435},
is essentially self-adjoint.

{}From the asymptotic behaviour we see that $wf(z)$ and $wg(z)$
are linearly independent solutions of the
recurrence in Lemma \ref{4180}, and moreover that
they both belong to $S^+_z$ for $q^2<c<1$.
The other statements follow easily from the asymptotics described
above.
\end{proof}
\end{sspar}

\begin{sspar}{\label{4190}} The Wronskian
$$
[wF(y), wF(y^{-1})] = \lim_{k\to-\infty}
a_k w_{k+1}w_k\bigl(
F_{k+1}(y)F_k(y^{-1})-F_k(y)F_{k+1}(y^{-1})\bigr)
=\hf (y^{-1}-y)
$$
using $a_k\to \hf$ as $k\to-\infty$ and the asymptotics
of \ref{4180}.
Note that the Wronskian is non-zero for $y\not= \pm 1$
or $z\not=\pm1$.
Since $wF(y)$ and $wF(y^{-1})$ are linearly independent
solutions to the
recurrence in Lemma \ref{4175}
for $z\in\C\backslash\R$, we see that we can express
$f_k(\mu(y))$ in terms of $F_k(y)$ and $F_k(y^{-1})$.
These solutions are related by the expansion
\begin{equation}
\begin{split}
f_k(\mu(y))
&= c(y) F_k(y) + c(y^{-1})F_k(y^{-1}), \\
c(y) &= \frac{(c/dy,d/y,dry,q/dry;q)_\infty}
{(y^{-2},c,r,q/r;q)_\infty},
\end{split}
\label{eq4191}
\end{equation}
for $c\not\in q^{-\Zp}$,
$y^2\not\in q^\Z$, see \cite[(4.3.2)]{GaspR}.
This shows that we have
$$
[wf(\mu(y)),wF(y)] =\hf c(y^{-1})(y-y^{-1}).
$$
\end{sspar}

\begin{sspar}{\label{4200}} Let us assume first that $0<c\leq q^2$,
so that $L$ is essentially self-adjoint. Then for
$z\in\C\backslash\R$ we have
$\phi_z = wf(z)$ and $\Phi_z=wF(y)$, where $z=\mu(y)$ and
$|y|<1$. In particular, it follows that $\phi_{x\pm i\ep}\to
\phi_x$ as $\ep\downarrow 0$.
For the asymptotic solution $\Phi_z$ we have
to be more careful in computing the limit of $z$ to the real axis.
For $x\in \R$ satisfying $|x|>1$ we have
$\Phi_{x\pm i\ep} \rightarrow wF_y$ as
$\ep\downarrow 0$, where
$y\in (-1,1)\backslash \{0\}$ is such that $\mu(y)=x$.
If $x\in [-1,1]$, then we put $x=\cos\chi=\mu(e^{i\chi})$
with $\chi\in [0,\pi]$, and then
$\Phi_{x-i\ep}\rightarrow wF_{e^{i\chi}}$
and $\Phi_{x+i\ep}\rightarrow wF_{e^{-i\chi}}$
as $\ep\downarrow 0$.
\end{sspar}

\begin{sspar}{\label{4210}} We calculate the integrand
in the Stieltjes-Perron inversion formula of \ref{2120}
using Proposition \ref{470} and (\ref{eq475}) in the
case $|x|<1$, where $x=\cos\chi=\mu(e^{i\chi})$.
For $u,v\in {\mathcal D}(\Z)$ we have
\begin{equation*}
\begin{split}
&\lim_{\ep\downarrow 0}\langle G(x+i\ep)u,v\rangle -
\langle G(x-i\ep)u,v\rangle = \\ & \lim_{\ep\downarrow 0}
\sum_{k\leq l}
\Bigl( \frac{(\Phi_{x+i\ep})_k(\phi_{x+i\ep})_l}
{[\phi_{x+i\ep},\Phi_{x+i\ep}]} -
\frac{(\Phi_{x-i\ep})_k(\phi_{x-i\ep})_l}
{[\phi_{x-i\ep},\Phi_{x-i\ep}]} \Bigr)
\bigl( u_l\bar v_k+u_k\bar v_l\bigr) (1-\hf
\de_{k,l})
= \\ &2\sum_{k\leq l}
\Bigl( \frac{w_kF_k(e^{-i\chi})w_lf_l(\cos\chi)}
{c(e^{i\chi})(e^{-i\chi}-e^{i\chi})} -
\frac{w_kF_k(e^{i\chi})w_lf_l(\cos\chi)}
{c(e^{-i\chi})(e^{i\chi}-e^{-i\chi})}\Bigr)
\bigl( u_l\bar v_k+u_k\bar v_l\bigr) (1-\hf
\de_{k,l})
= \\ &2\sum_{k\leq l}
\Bigl( w_kw_lf_l(\cos\chi)
\frac{c(e^{-i\chi})F_k(e^{-i\chi})+c(e^{i\chi})F_k(e^{i\chi})}
{c(e^{i\chi})c(e^{-i\chi})(e^{-i\chi}-e^{i\chi})}
\bigl( u_l\bar v_k+u_k\bar v_l\bigr) (1-\hf
\de_{k,l}) = \\
&2\sum_{k\leq l}
\Bigl( \frac{w_kw_lf_l(\cos\chi)f_k(\cos\chi)}
{c(e^{i\chi})c(e^{-i\chi})(e^{-i\chi}-e^{i\chi})}
\bigl( u_l\bar v_k+u_k\bar v_l\bigr) (1-\hf
\de_{k,l}) = \\
&\frac{2}{c(e^{i\chi})c(e^{-i\chi})(e^{-i\chi}-e^{i\chi})}
\sum_{l=-\infty}^\infty w_lf_l(\cos\chi)u_l
\sum_{k=-\infty}^\infty w_kf_k(\cos\chi)\bar v_k
\end{split}
\end{equation*}
using the expansion (\ref{eq4191}) and the Wronskian in
\ref{4190}. Now integrate over the interval
$(a,b)$ with $-1\leq a<b\leq 1$ and replacing $x$ by $\cos\chi$, so
that $\frac{1}{2\pi i}dx = (e^{i\chi}-e^{-i\chi})d\chi/4\pi$
we obtain, with $a=\cos \chi_a$, $b=\cos\chi_b$, and
$0\leq \chi_b<\chi_a\leq\pi$,
\begin{equation*}
\begin{split}
E_{u,v}\bigl( (a,b)\bigr) &= \frac{1}{2\pi}
\int_{\chi_b}^{\chi_a} \bigl({\mathcal F}u\bigr)(\cos\chi)
\overline{\bigl({\mathcal F}v\bigr)(\cos\chi)}
\frac{d\chi}{|c(e^{i\chi})|^2}, \\
\bigl({\mathcal F}u\bigr)(x) &= \langle u, \phi_x\rangle =
\sum_{l=-\infty}^\infty w_lf_l(\cos\chi)u_l.
\end{split}
\end{equation*}
This shows that $[-1,1]$ is contained in the continuous
spectrum of $L$.
\end{sspar}

\begin{sspar}{\label{4220}} For $|x|>1$ we can calculate as
in \ref{4210} the integrand in the Stieltjes-Perron
inversion formula, but now we use that $x=\mu(y)$ with
$|y|<1$, see \ref{4200}. This gives
$$
\lim_{\ep\downarrow 0}\langle G(x+i\ep)u,v\rangle
= 2\sum_{k\leq l} \frac{w_kF_k(y)w_lf_l(y)}
{c(y^{-1})(y-y^{-1})}
\bigl( u_l\bar v_k+u_k\bar v_l\bigr) (1-\hf
\de_{k,l}),
$$
and since the limit
$\lim_{\ep\downarrow 0}\langle G(x+i\ep)u,v\rangle$ gives the
same result, we see, as in the case of the Meixner
functions, that we can only have discrete mass points
for $|x|>1$ in the spectral measure at the zeroes of the
Wronskian, i.e. at the zeroes of $y\mapsto c(y^{-1})$
with $|y|<1$ or at $y=\pm 1$. Let us assume that all zeroes
of the $c$-function are simple, so that the spectral measure
at these points can be easily calculated.

The zeroes of the $c$-function can be read off from the
expressions in (\ref{eq4191}), and they are
$\{ cq^k/d\mid k\in\Zp\}$, $\{ dq^k\mid k\in\Zp\}$ and
$\{ q^k/dr\mid k\in\Z\}$. Assuming that $|c/d|<1$ and
$|d|<1$, we see that the first two sets do not contribute.
In the more general case we have that the product is
less than $1$, since the product equals $c<1$. We leave this
extra case to the reader.
The last set, labeled by $Z$ always contributes to
the spectral measure. Now
for $u,v\in{\mathcal D}(\Z)$ we let $x_p=\mu(y_p)$,
$y_p=q^p/dr$, $p\in\Z$, with
$|q^p/dr|>1$, so that, cf. \ref{4140},
$$
E_{u,v}(\{ x_p\}) = \text{Res}_{y=y_p^{-1}}
\Bigl( \frac{-1}{c(y^{-1})y}\Bigr)
w_k F_k(y_p^{-1}) w_l f_l(x_p)
\bigl( u_l\bar v_k+u_k\bar v_l\bigr) (1-\hf
\de_{k,l})
$$
after substituting $x=\mu(y)$. Now from (\ref{eq4191})
we find $f_k(x_p)=c(y_p)F_k(y_p^{-1})$, since $c(y_p^{-1})=0$
and we assume here that $c(y_p)\not=0$. Hence, we
can symmetrise the sum again and find
$$
E_{u,v}(\{ x_p\}) = \Bigl(\text{Res}_{y=y_p}
\frac{1}{c(y^{-1})c(y)y}\Bigr)
\bigl({\mathcal F}u\bigr)(x_p) \overline{\bigl({\mathcal F}v\bigr)(x_p)}
$$
switching to the residue at $y_p$.
\end{sspar}

\begin{sspar}{\label{4230}} We can combine the calculations
in the following theorem. Note that most of the regularity
conditions can be removed by continuity after calculating
explicitly all the residues. The case of an extra set of
finite mass points is left to the reader, cf.
\ref{4220}, as well as the
case of other choices of the parameters $c$, $d$ and $r$
for which the expression under the square root sign in $a_k$
in \ref{4170} is positive. See \cite[App.~A]{KoelSsu}
for details.

\begin{thm} Assume $r<0$, $0<c\leq q^2$, $d\in\R$
with $|d|<1$ and $|c/d|<1$ such that
the zeroes of $y\mapsto c(y)$ are simple and $c(y)=0$
implies $c(y^{-1})\not= 0$. Then the spectral measure
for the Jacobi operator on $\Hi$ defined by {\rm \ref{4170}}
is given by, $A\subset\R$ a Borel set,
\begin{equation*}
\begin{split}
&\langle E(A)u,v\rangle =
\int_{\cos\chi\in [-1,1]\cap A}\bigl({\mathcal F}u\bigr)(\cos\chi)
\overline{\bigl({\mathcal F}v\bigr)(\cos\chi)} \frac{d\chi}
{|c(e^{i\chi})|^2} \\ &+
\sum_{p\in\Z, |q^p/dr|>1, \mu(q^p/dr)\in A}
\Bigl(\text{Res}_{y=q^p/dr}
\frac{1}{c(y^{-1})c(y)y}\Bigr)
\bigl({\mathcal F}u\bigr)(\mu(q^p/dr))
\overline{\bigl({\mathcal F}v\bigr)(\mu(q^p/dr))}.
\end{split}
\end{equation*}
\end{thm}

\begin{proof} It only remains to prove that $\pm 1$ is not
contained in the point spectrum. These are precisely the points
for which $F(y)$ and $F(y^{-1})$ are not linearly independent
solutions. We have to show that $\phi_{\pm 1}\not\in\Hi$, and
this can be done by determining its asymptotic behaviour
as $k\to-\infty$, see \cite{Kake}, \cite{KoelSbig} for more
information.
\end{proof}

Take $A=\R$ and $u=e_k$ and $v=e_l$, then we find the following
orthogonality relations for the ${}_2\vp_1$-series as
in Lemma \ref{4175};
\begin{equation*}
\begin{split}
&\int_0^\pi f_k(\cos\chi) f_l(\cos\chi)
\frac{d\chi}{|c(e^{i\chi})|^2} +
\\ &\sum_{p\in\Z, |q^p/dr|>1}
\Bigl({\text{Res}}_{y=q^p/dr}
\frac{1}{c(y^{-1})c(y)y}\Bigr)
f_k(\mu(\frac{q^p}{dr}))f_l(\mu(\frac{q^p}{dr}))
= \frac{\de_{k,l}}{w_k^2}.
\end{split}
\end{equation*}
\end{sspar}

\begin{sspar}{\label{4240}} In Theorem \ref{4230} we have
made the assumption $0<c\leq q^2$ in order to have $L$
as an essentially self-adjoint operator. From the general
considerations in \ref{450}-\ref{472} it follows
that the previous calculations, and in particular
Theorem \ref{4230}, remain valid for $q^2<c<1$ if we can show
that there exists a self-adjoint extension of $L$ satisfying
the assumptions of \ref{460}. It suffices to check (2) of
\ref{460}, or, by Lemma \ref{440},
the existence of a $\te\in[0,2\pi)$ such that
\begin{equation}
\lim_{k\to\infty} [wf(z), e^{i\te}wF(i(1-\sqrt{2}))
+ e^{-i\te}wF(-i(1-\sqrt{2}))]_k = 0.
\label{eq4241}
\end{equation}
Indeed, $wF(\pm i(1-\sqrt{2}))$ is the element $\Phi_{\pm i}$
up to the normalisation of the length. This is not important for
showing the existence of $\te$ satisfying (\ref{eq4241}).

For (\ref{eq4241}) we need to know the asymptotic behaviour
of $F_k(y)$ as $k\to\infty$. The same result in basic
hypergeometric series that results in (\ref{eq4191}) can be used to
prove that
\begin{equation*}
\begin{split}
&\qquad\qquad F_k(y) = a(y)\, f_k(y) + b(y)\, g_k(y), \\
&a(y) = \frac{(qdy/c,qy/d,qcy/dr,dr/cy;q)_\infty}
{(qy^2,q/c,qc/d^2r,d^2r/c;q)_\infty}, \quad
b(y) = \frac{(dy,cy/d,q^2y/dr,dr/yq;q)_\infty}
{(qy^2,c/q,qc/d^2r,d^2r/c;q)_\infty},
\end{split}
\end{equation*}
for $y^2,c\notin q^\Z$,
so that $F_k(y)=a(y){\mathcal O}(1)+b(y){\mathcal O}((q/c)^k)$.
It follows that
$$
\lim_{k\to\infty} [wf(z),wF(y)]_k =
\hf \big\vert \frac{c}{dr}\big\vert (1-\frac{c}{q})\, b(y),
$$
so that the limit in (\ref{eq4241}) equals
$$
\frac{\big\vert \frac{c}{dr}\big\vert (1-\frac{c}{q})}
{2 \| F(i(1-\sqrt{2}))\|}
\bigl( e^{i\te}b(i(1-\sqrt{2})) +e^{-i\te}b(-i(1-\sqrt{2}))
\bigr).
$$
The term in parentheses is $2\Re(e^{i\te}b(i(1-\sqrt{2})))$,
which is zero for $\te = -\frac{\pi}{2} + \arg b(i(1-\sqrt{2}))$.
\end{sspar}

\begin{sspar}{\label{4250}} The Meixner functions can
formally be obtained as a limit case as $q\uparrow 1$ from
the spectral analysis of the second order $q$-difference
operator as considered here. For this we make the
following specialisation; $c\mapsto qs^{-2}$,
$d\mapsto q^{1+\la}s^{-1}$, $r\mapsto q^{-\ep-\la}$, where
$\ep$ and $\la$ have the same meaning as in
\ref{480}. Note that $r<0$ is no longer valid but
the operator is well-defined under suitable conditions
on $\la$, cf. \ref{480}. The
operator $L$ with \ref{4170} has the same type of spectral
measure (exercise). Now consider the operator
$$
L_q = \frac{2L-s-s^{-1}}{1-q}, \qquad
L_q f_k = a^q_k\, f_{k+1} + b^q_k\, f_k + a^q_{k-1}\, f_{k-1}
$$
for $f_k=e_{-k}$, then a calculation shows that
$a^q_k\to a_k$, $b^q_k\to b_k$ with $a_k$ and $b_k$ as in
\ref{480} with $2a=s+s^{-1}$ as $q\uparrow 1$. We now
assume $a\geq 1$ as in the previous subsection on
Meixner functions.

Now the operator $L_q$
has continuous spectrum supported on
$[(-2-s-s^{-1})/(1-q),(2-s-s^{-1})/(1-q)]$. For
$s+s^{-1}>2$ the continuous spectrum will disappear to $-\infty$
as $q\uparrow 1$; it will tend to $(-\infty,0]$ for $s+s^{-1}=2$,
and it will tend to $\R$ if $0\leq s+s^{-1}<2$.
The discrete spectrum (or at least the infinite number of
discrete mass points) is of the form
$$
s\frac{q^{p-1+\ep}-1}{1-q} + s^{-1} \frac{q^{1-\ep-p}-1}{1-q},
\qquad p\in\Z,\ |sq^{p+1-\ep}|>1.
$$
As $q\uparrow 1$ this tends to $(p+\ep-1)(s^{-1}-s)$
with $p\in\Z$ for $|s|>1$ and it will disappear for
$|s|\leq 1$. So for $a>1$ only the discrete spectrum survives
and this corresponds precisely to Theorem \ref{4140}.
For $0\leq a\leq 1$ the discrete spectrum disappears in the
limit.

It is also possible to show that the solutions
for the $q$-hypergeometric difference operator tend to
to solutions of the operator for the Meixner functions.
This requires the transformation as in \ref{4115},
and one of Heine's transformations.

So we have motivated, at least formally, that the
orthogonality relations for the
Meixner functions
can be obtained as a limit case of the ${}_2\vp_1$-series that
arise as solutions of the second order $q$-difference
hypergeometric equation, and that the limit
transition remains valid on the level of (the support of)
the spectral measure. For the Laguerre and 
Meixner-Pollaczek
functions the limit remains valid on the level of the
support of the spectral measure.
\end{sspar}


\newpage
\bibliographystyle{amsplain}

\end{document}